Article

# Structure of Finite-Dimensional Protori

Wayne Lewis

University of Hawaii, 874 Dillingham Blvd., Honolulu, HI 96817, USA; waynel@math.hawaii.edu

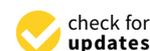

**Abstract:** A *Structure Theorem for Protori* is derived for the category of finite-dimensional *protori* (compact connected abelian groups), which details the interplay between the properties of density, discreteness, torsion, and divisibility within a finite-dimensional protorus. The spectrum of resolutions for a finite-dimensional protorus are parameterized in the structure theorem by the dual category of finite rank torsion-free abelian groups. A consequence is a *universal resolution* for a finite-dimensional protorus, independent of a choice of a particular subgroup. A resolution is also given strictly in terms of the path component of the identity and the union of all zero-dimensional subgroups. The structure theorem is applied to show that a morphism of finite-dimensional protori lifts to a product morphism between products of periodic locally compact groups and real vector spaces.

## 1. Introduction

All compact groups herein are finite-dimensional and all torsion-free groups have finite rank. We carry out a study of the structure of compact, connected abelian groups, or *protori*. A practical new perspective on the category of torsion-free abelian groups, from the point of view of the dual category of protori, emerges organically from a detailed analysis of the algebro-topological structure of protori. A Structure Theorem for Protori (Theorem 1) is derived which applies to all objects in the category of protori.

Rather than a study involving specialization and classification relative to a particular morphism, we take a holistic approach to the category of protori. The main results are a Structure Theorem for Protori (Theorem 1), a universal resolution for a protorus (Corollary 6), a structural result on the lattice of compact open subgroups of zero-dimensional subgroups of a protorus under a natural locally compact topology (Proposition 6), and a lifting theorem for morphisms of protori (Theorem 2), which facilitates a reduction to a decoupled analysis of morphisms between periodic LCA groups.

The Structure Theorem for Protori is formulated for an arbitrary protorus by applying the key new Lemma 5, intrinsically engineered for protori, to the *Resolution Theorem for Compact Abelian Groups* (Proposition 2.2, [1]), which states that a compact abelian group $H$ is topologically isomorphic to $[\Delta \times \mathfrak{L}(H)]/\Gamma$ for a totally disconnected $\Gamma$ and a profinite subgroup $\Delta$ of $H$ such that $H/\Delta$ is a torus, where $\mathfrak{L}(H)$ is the *Lie algebra* of $H$ (Proposition 7.36, [2]). An immediate consequence of the Structure Theorem for Protori is the existence of a *universal resolution* for a protorus $G$ (Corollary 6): $[\widehat{\Delta}_{X_\infty} \times \mathfrak{L}(G)]/X_\infty$, where $\widehat{\Delta}_{X_\infty}$ is a *periodic* (Definition 1.13, [3]), locally compact *topological divisible hull* of a profinite $\Delta$ of a given resolution of $G$, $X_\infty$ is a minimal quotient-divisible extension of an intervening Pontryagin dual of $G$, and the concept of *minimal divisible locally compact cover* of $G$ is introduced (Corollary 7) and realized via $\widehat{\Delta}_{X_\infty} \times \mathfrak{L}(G)$. The canonical zero-dimensional periodic group $\widehat{\Delta}_{X_\infty}$ is an inverse limit of discrete groups, topologically isomorphic to the local direct product of its divisible *p*-Sylow components, each isomorphic to a product of powers of the *p*-adic numbers $\widehat{\mathbb{Q}}_p$ and/or the Prüfer group $\mathbb{Z}(p^\infty)$ (pp. 48–49, [3]). Proposition 4 details a fundamental new result deconstructing periodic LCA groups intrinsic to protori; the result facilitates an approach enabling one to assess the impact in situ of addition to a fixed profinite $\Delta^*$ a torsion-free $Y^*$ with rk $Y^*$ = dim $G$, effecting a parametrization of topologically isomorphic resolutions via the spectrum of such torsion-free groups.

We would be remiss not to toot the horn a bit, emphasizing Theorem 1 gives for the first time resolutions of a protorus *G* not in terms of a particular profinite subgroup in L(*G*) but rather in terms of its topological divisible hull in *G*; further, a topological divisible hull of some $\Delta \in $ L(*G*) serves merely as an upper bound for the spectrum of resolutions associated with torsion-free subgroups spanning the region between a dense free subgroup $\mathbb{Z}_\Delta$ of $\Delta$ and a minimal quotient-divisible extension of the Pontryagin dual of *G*; in the large, the Structure Theorem for Protori gives a pseudoalgorithm for deriving a periodic LCA group and associated resolution determined by *any* torsion-free group bounded below by $\mathbb{Z}_\Delta$, $\Delta \in $ L(*G*), and above by the divisible hull of $\mathbb{Z}_\Delta$ in *G*. The Structure Theorem for Protori not only describes the structure of an arbitrary protorus, but it provides a tool, with instructions, for surgically deconstructing protori.

Applying the structure theory developed in Theorem 1, we derive a Structure Theorem for Morphisms (Theorem 2) a new result stating that a morphism $f: G \to H$ of protori with duals *X* and *Y* lifts to a product map between minimal divisible locally compact covers $f \mid_{\widehat{\Delta}} \times f_{\mathfrak{L}}: \widehat{\Delta}_{X_\infty} \times \mathfrak{L}(G) \to \widehat{\Delta}_{Y_\infty} \times \mathfrak{L}(H)$. Because $\mathfrak{L}(G)$ and $\mathfrak{L}(H)$ are finite-dimensional real topological vector spaces (Proposition 7.24, [3]), the Structure Theorem for Morphisms reduces the analysis of protori morphisms to those between divisible periodic LCA groups. Because $\widehat{\Delta}_{X_\infty}$ and $\widehat{\Delta}_{Y_\infty}$ are divisible periodic LCA groups topologically isomorphic to the local product of their *p*-Sylow components, respected by protori morphisms, protori morphisms are an amalgam of their restrictions to *p*-Sylow factors, where the action on such a factor is uniquely determined by its action on a compact open subgroup. Thus, the Structure Theorem for Morphisms implies that deconstructing protori morphisms effectively reduces to the analysis of morphisms between finitely generated $\widehat{\mathbb{Z}}_p$-modules. Lastly, we state without proof that Theorem 2 generalizes the analogous result that a morphism between complex tori lifts to a complex-linear map between complex vector spaces (Proposition 2.1, [4]).

Regarding the breakdown of sections comprising the paper: Section 2 provides the requisite background for our study. Section 3 proves the main structure theory results. Section 4 gives several illustrative applications involving *projective resolutions*, *ACD groups*, and *morphisms of protori*, culminating in a general lifting theorem for the category of protori.

## 2. Background

A ***protorus*** is a compact connected abelian group. The name *protorus* derives from the formulation of its definition as an inverse limit of tori (Corollary 8.18, Proposition 1.33) [2], analogous to a *profinite* group as an inverse limit of finite groups. A ***morphism*** between topological groups is a continuous homomorphism. A ***topological isomorphism*** is an open bijective morphism between topological groups, which we indicate by $\cong_t$. Set $\mathbb{T} \stackrel{\text{def}}{=} \mathbb{R}/\mathbb{Z}$ with the quotient topology induced from the Euclidean topology on $\mathbb{R}$, for which $\mathbb{Z}$ is discrete under the subspace topology. A ***torus*** is a topological group topologically isomorphic to $\mathbb{T}^n$ for some positive integer *n*. A protorus is ***torus-free*** if it contains no subgroups topologically isomorphic to a torus.

All groups herein are abelian and all topological groups are Hausdorff. All finite-dimensional real topological vector spaces are topologically isomorphic to a real Euclidean vector space of the same dimension. All references to duality refer to Pontryagin duality. Finitely generated in the context of profinite groups will always mean topologically finitely generated. If *A* and *B* are topological groups which each contain an isomorphic copy $X^*$ of a torsion-free group *X* such that $X^*$ embeds diagonally in $A \times B$ as a closed subgroup, then we write $(A \times B)/X$ for the associated quotient. Some authors use the term *solenoid* or *solenoidal group* to describe protori; we use *protorus* to connote compact connected abelian group and *solenoid* to mean one-dimensional protorus.

Pontryagin duality is a contravariant endofunctor on the category of locally compact abelian groups under continuous homomorphism, $^\vee : \mathbb{LCA} \to \mathbb{LCA}$, given by $G^\vee = \underset{\text{continuous}}{\mathrm{Hom}}(G, \mathbb{T})$ under the topology of compact convergence and $\rho^\vee : H^\vee \to G^\vee$ by $\rho^\vee(\chi) = \chi \circ \rho$ for $\rho : G \to H$, such that $^{\vee\vee}$ is naturally isomorphic to the identify functor. Each object in the category is isomorphic to some image of the Pontryagin duality functor. Pontryagin duality restricts to an equivalence between the category of discrete abelian groups and the opposite category of compact abelian groups (Theorem 7.63) [2] wherein compact abelian groups are connected if and only if they are divisible (Proposition 7.5(i), [2]), (24.3, [5]), (23.17, [5]), (24.25, [5]). Some locally compact abelian groups, such as finite cyclic groups $\mathbb{Z}(n)$, the real numbers $\mathbb{R}$, the *p*-adic numbers $\hat{\mathbb{Q}}_p$, and the adeles $\mathbb{A}$ are categorical fixed points of the Pontryagin duality functor.

For a compact abelian group $G$, the ***Lie algebra*** $\mathfrak{L}(G) \overset{\text{def}}{=} \underset{\text{continuous}}{\mathrm{Hom}}(\mathbb{R}, G)$, consisting of the set of continuous homomorphisms under the topology of compact convergence, is a real topological vector space (Proposition 7.36, [2]). The ***exponential function*** of $G$, $\exp_G : \mathfrak{L}(G) \to G$ given by $\exp_G(r) = r(1)$, is a morphism of topological groups, and $\exp_G$ is injective when $G$ is torsion-free (Corollary 8.47, [2]). Let $G_0$ denote the ***connected component of the identity*** and $G_a = \exp_G \mathfrak{L}(G)$ the ***path component of the identity*** in $G$ (Theorem 8.30, [2]).

The ***dimension*** of a compact abelian group $G$ is $\dim G \overset{\text{def}}{=} \dim_\mathbb{R} \mathfrak{L}(G)$. When $G$ is a compact abelian group, $\mathfrak{L}(G) \cong_t \mathbb{R}^{\dim G}$ as topological vector spaces (Proposition 7.24, [2]) and $\dim G = \dim_\mathbb{Q}(\mathbb{Q} \otimes_\mathbb{Z} G^\vee)$ when $G$ has positive dimension (Theorem 8.22, [2]). A sequence of compact abelian groups $G_1 \overset{\phi}{\rightarrowtail} G_2 \overset{\psi}{\twoheadrightarrow} G_3$ is ***exact*** if $\phi$ and $\psi$ are, respectively, injective and surjective morphisms and $\mathrm{Ker}\, \psi = \mathrm{Im}\, \phi$; note that automatically $\phi$ is open onto its image and $\psi$ is open (Theorem 5.29, [5]). For a morphism $\rho : G \to H$ of locally compact abelian groups, the ***adjoint*** of $\rho$ is the morphism $\rho^\vee : H^\vee \to G^\vee$ given by $\rho^\vee(\chi) = \chi \circ \rho$ (Theorem 24.38, [5]). A sequence of compact abelian groups $G_1 \overset{\phi}{\rightarrowtail} G_2 \overset{\psi}{\twoheadrightarrow} G_3$ is exact if and only if $G_3^\vee \overset{\psi^\vee}{\rightarrowtail} G_2^\vee \overset{\phi^\vee}{\twoheadrightarrow} G_3^\vee$ is an exact sequence of discrete abelian groups (Theorem 2.1, [6]). A compact abelian group $G$ is totally disconnected $\Leftrightarrow \dim G = 0 \Leftrightarrow G^\vee$ is torsion $\Leftrightarrow \dim(\mathbb{Q} \otimes G^\vee) = 0$ (Corollary 8.5, [2]).

Torsion-free abelian groups $A$ and $B$ are ***quasi-isomorphic*** if there is $f : A \to B$, $g : B \to A$, and $0 \neq n \in \mathbb{Z}$ such that $fg = n \cdot 1_B$ and $gf = n \cdot 1_A$; $A$ and $B$ are ***nearly-isomorphic*** if for each $0 \neq n \in \mathbb{Z}$ there is a relatively prime $m \in \mathbb{Z}$, $f : A \to B$, and $g : B \to A$ such that $fg = m \cdot 1_B$ and $gf = m \cdot 1_A$. By (Corollary 7.7, [7]), $A$ and $B$ are quasi-isomorphic if and only if there is a monomorphism $h : A \to B$ such that $A/f(B)$ is finite. It follows by Pontryagin duality that $A$ and $B$ are quasi-isomorphic if and only if there is a surjective morphism $h^\vee : B^\vee \to A^\vee$ with finite kernel. This is exactly the definition of ***isogeny*** between protori: $G$ and $H$ are ***isogenous*** if there is a surjective morphism $G \to H$ with finite kernel. Define protori $G$ and $H$ to be ***topologically nearly-isomorphic*** if $G^\vee$ and $H^\vee$ are nearly isomorphic. Evidently from the definitions, *quasi-isomorphism and near-isomorphism of torsion-free groups* and *isogeny and topological near-isomorphism of protori* are equivalence relations.

Let $\mathbb{P}$ denote the set of prime numbers. A ***supernatural number*** is a formal product $\mathbf{n} = \prod_{p \in \mathbb{P}} p^{\mathbf{n}_p}$ where $0 \leq \mathbf{n}_p \in \mathbb{Z} \cup \{\infty\}$ (Section 2.3, [8]). We denote the *p*-adic integers $\hat{\mathbb{Z}}_p$ and the ***p*-adic numbers** $\hat{\mathbb{Q}}_p$. Let $\mathbb{S}$ denote the set of all supernatural numbers. A ***height sequence*** $(\mathbf{n}_p)_{p \in \mathbb{P}}$ is a sequence of exponents of a supernatural number. Define an equivalence relation on $\mathbb{S}$ by stipulating that supernatural numbers $\mathbf{m}$ and $\mathbf{n}$ are ***equivalent*** if their height sequences are equal except for a finite number of primes $p$ for which $\mathbf{m}_p, \mathbf{n}_p < \infty$. Define a ***type*** to be an equivalence class of a height sequence, denoted $\tau = [(\mathbf{n}_p)_{p \in \mathbb{P}}]$. For an element $a$ of a torsion-free group $X$, the ***p-height of a in X***, $\mathrm{ht}_p^X(a)$, $p \in \mathbb{P}$, is $n$ if there exists $0 \leq n \in \mathbb{Z}$ such that $a \in p^n X \backslash p^{n+1} X$ and $\infty$ otherwise. The ***height sequence of a in X*** is $(\mathrm{ht}_p^X(a))_{p \in \mathbb{P}}$. Define the ***type of a in X*** to be $[(\mathrm{ht}_p^X(a))_{p \in \mathbb{P}}]$. Because any two nonzero elements of a rank-1 torsion-free group $A$ have the same type in $A$, the ***type*** of $A$, $\mathrm{type}(A)$, is the type of any nonzero element in $A$. Two rank-1 groups $A$ and $B$ are isomorphic if and only if $\text{type}(A) = \text{type}(B)$, and given $\mathbf{n} \in \mathbb{S}$ there is a rank-1 torsion-free group $C$ with $\text{type}(C) = [(n_p)_{p \in \mathbb{P}}]$ (Theorem 1.1, [7]). For a torsion-free group $X$, set $\tau_{\sup}(X) \stackrel{\text{def}}{=} [(\sup\{\text{ht}_p^X(a) : 0 \neq a \in X\})_{p \in \mathbb{P}}]$.

An abelian group $D$ is **divisible** if for every $0 \neq d \in D$ and $0 \neq k \in \mathbb{Z}$ there is a $d' \in D$ such that $kd' = d$. A torsion-free group $X$ is **quotient-divisible** if it contains a free subgroup $F$ such that $X/F$ is a divisible (torsion) group (p. 473, [9]). Define the **radius** of $0 \neq z = (z_1, \ldots, z_n) \in \mathbb{Z}^n$, $0 < n \in \mathbb{Z}$, to be $\mathbf{r}(z) \stackrel{\text{def}}{=} \gcd(z_1, \ldots, z_n)$ and set $\mathbf{r}(0) \stackrel{\text{def}}{=} 0$. Set $\mathbb{B}^n(r) \stackrel{\text{def}}{=} \{z \in \mathbb{Z}^n : \mathbf{r}(z) \leq r\}$ and define a **unit hemisphere** $H \subset \mathbb{Z}^n$ to be a subset of $\mathbb{B}^n(1) \setminus \{0\}$ for which each line through 0 in $\mathbb{Q}^n$ passes through exactly one point in $H$. A minimal quotient-divisible extension of a rank-$n$ torsion-free group $X$ with $H \subseteq \mathbb{Z}^n \subseteq X \subseteq \mathbb{Q}^n$, $H$ a unit hemisphere, is $X_\infty \stackrel{\text{def}}{=} \sum_{z \in H} \sum_{p \in \mathbb{P}} \sum_{\text{ht}_p^X(z) > 0} \mathbb{Z}[\frac{1}{p}]z$ where $\mathbb{Z}[\frac{1}{p}] = \{\pm \frac{j}{p^\ell} : 0 \leq j, \ell \in \mathbb{Z}\}$.

A **profinite group** is an inverse limit of finite groups or, equivalently, a totally disconnected compact Hausdorff group (Theorem 1.34, [2]); such a group is **finitely generated** if it contains a dense finitely generated subgroup. Profinite abelian groups $D$ and $E$ are **isogenous** if there are morphisms $f: D \to E$ and $g: E \to D$ such that $E/f(D)$ and $D/g(E)$ are bounded torsion groups; for finitely generated $D$ and $E$, this is equivalent to $E/f(D)$ and $D/g(E)$ being finite. By symmetry, isogeny of profinite abelian groups is an equivalence relation. Proceeding strictly according to Pontryagin duality, one would conclude that torsion abelian groups $A$ and $B$ be defined as **quasi-isomorphic** if there are morphisms $h: A \to B$ and $k: B \to A$ such that $B/h(A)$ and $A/k(B)$ are bounded torsion groups; this is, in fact, the definition for quasi-isomorphism between torsion abelian groups (see Proposition 1.8, [7]).

The development of a structure theory for protori is very much dependent on the theory of profinite abelian groups. The profinite theory in this section is derived in large part from the standard reference by Ribes and Zaleeskii [8]. We begin by showing that the additivity of dimension for vector spaces also holds for compact abelian groups.

**Lemma 1.** *If $0 \to G_1 \to G_2 \to G_3 \to 0$ is an exact sequence of compact abelian groups, then $\dim G_2 = \dim G_1 + \dim G_3$.*

**Proof.** The exactness of $0 \to G_1 \to G_2 \to G_3 \to 0$ implies the exactness of $0 \to G_3^\vee \to G_2^\vee \to G_1^\vee \to 0$ and this implies the exactness of $0 \to \mathbb{Q} \otimes_\mathbb{Z} G_3^\vee \to \mathbb{Q} \otimes_\mathbb{Z} G_2^\vee \to \mathbb{Q} \otimes_\mathbb{Z} G_1^\vee \to 0$ because $\mathbb{Q}$ is torsion-free (Theorem 8.3.5, [9]). However, this is an exact sequence of $\mathbb{Q}$-vector spaces and hence $\dim_\mathbb{Q}(\mathbb{Q} \otimes_\mathbb{Z} G_2) = \dim_\mathbb{Q}(\mathbb{Q} \otimes_\mathbb{Z} G_3) + \dim_\mathbb{Q}(\mathbb{Q} \otimes_\mathbb{Z} G_1)$. This establishes the claim because, in general $\dim G = \dim_\mathbb{Q} \mathbb{Q} \otimes_\mathbb{Z} G^\vee$ by (Theorem 8.22, [2]) for $\dim G \geq 1$ and $\dim G = 0 \Leftrightarrow \dim_\mathbb{Q}(\mathbb{Q} \otimes_\mathbb{Z} G^\vee) = 0$. □

Fix $n \in \mathbb{Z}$. Denote by $\mu_n$ the *multiplication-by-$n$* map $A \to A$ for an abelian group $A$, given by $\mu_n(a) = na$ for $a \in A$.

**Lemma 2.** $\mu_n : G \to G$, $0 \neq n \in \mathbb{Z}$, *is an isogeny for a protorus $G$.*

**Proof.** $\mu_n$ is a surjective morphism because $G$ is a divisible abelian topological group, so the adjoint $\mu_n^\vee : G^\vee \to G^\vee$ is injective, whence $[G^\vee : \mu_n^\vee(G^\vee)]$ is finite by (Proposition 6.1.(a), [7]). It follows that $\ker \mu_n$ is finite and $\mu_n$ is an isogeny. □

A profinite group is either finite or uncountable (Proposition 2.3.1, [8]). The **profinite integers** $\widehat{\mathbb{Z}}$ is the inverse limit of cyclic groups of order $n$. $\widehat{\mathbb{Z}}$ is topologically isomorphic to $\prod_{p \in \mathbb{P}} \widehat{\mathbb{Z}}_p$ (Example 2.3.11, [8]) and to the profinite completion of $\mathbb{Z}$ (Example 2.1.6.(2), [8]), whence $\widehat{\mathbb{Z}}^m$ is a finitely generated profinite abelian group, $0 \leq m \in \mathbb{Z}$.

For a protorus $G$, the *Resolution Theorem for Compact Abelian Groups* states that $G$ contains a profinite subgroup $\Delta$ such that $G \cong_t \frac{\Delta \times \mathfrak{L}(G)}{\Gamma}$ where $\Gamma$ is a discrete subgroup of $\Delta \times \mathfrak{L}(G)$ and $G/\Delta \cong_t \mathbb{T}^{\dim G}$ (Theorems 8.20 and 8.22, [2]). In this case, the exact sequence $\Delta \rightarrowtail G \twoheadrightarrow \mathbb{T}^{\dim G}$ dualizes to $\mathbb{Z}^{\dim G} \hookrightarrow G^{\vee} \twoheadrightarrow \Delta^{\vee}$ where, without loss of generality, $G^{\vee} \subseteq \mathbb{Q}^{\dim G}$ so that $\Delta^{\vee} \cong \frac{G^{\vee}}{\mathbb{Z}^{\dim G}} \subseteq \frac{\mathbb{Q}^{\dim G}}{\mathbb{Z}^{\dim G}} \cong (\frac{\mathbb{Q}}{\mathbb{Z}})^{\dim G}$, whence by duality there is an epimorphism $\widehat{\mathbb{Z}}^{\dim G} \twoheadrightarrow \Delta$, because $\widehat{\mathbb{Z}} \cong_t (\frac{\mathbb{Q}}{\mathbb{Z}})^{\vee}$ (Example 2.9.5, [8]). It follows that, in the context of protori, the profinite groups of the *Resolution Theorem* are simultaneously finitely generated profinite abelian groups and finitely generated profinite $\widehat{\mathbb{Z}}$-modules. The continuous scalar multiplication $\mathbb{Z} \times \Delta \to \Delta$ is componentwise: if $\mathbf{x} = (x_1, \ldots, x_m) \in \Delta$, where $\mathbf{x}_j = (x_{j_p})_{p \in \mathbb{P}}$, then $k\mathbf{x} = (kx_1, \ldots, kx_m)$ where the scalar multiplication in each coordinate is given by $k\mathbf{x}_j = (kx_{j_p})_{p \in \mathbb{P}}$, applying the usual scalar multiplications for $\widehat{\mathbb{Z}}_p$ and $\mathbb{Z}(p^r)$. A locally compact abelian group $K$ for which $\overline{\mathbb{Z}g}$ is profinite for each $x \in K$ contains a unique *p-Sylow subgroup*, $p \in \mathbb{P}$ (Theorem 3.3, [3]), and a profinite group $H$ can be decomposed uniquely into the product of its *p*-Sylow subgroups (Proposition 2.3.8, [8]).

**Lemma 3.** *The algebraic structure of a finitely generated profinite abelian group uniquely determines its topological structure.*

**Proof.** A profinite group has a neighborhood basis at 0 consisting of open (whence closed) subgroups (Theorem 1.34, [2]). A subgroup of a finitely generated profinite abelian group is open if and only if it has finite index (Lemma 2.1.2, Proposition 4.2.5, [8]). It follows that finitely generated profinite abelian groups are topologically isomorphic if and only if they are isomorphic as abelian groups. □

In light of Lemma 3, we usually write $\cong$ in place of $\cong_t$ when working with profinite subgroups of protori.

Set $\mathbb{Z}(p^n) \stackrel{\text{def}}{=} \frac{\mathbb{Z}}{p^n\mathbb{Z}}$ for $0 \leq r \in \mathbb{Z}$. We introduce the notation $\widehat{\mathbb{Z}}(p^n) \stackrel{\text{def}}{=} \frac{\widehat{\mathbb{Z}}_p}{p^n\widehat{\mathbb{Z}}_p} \cong \frac{\mathbb{Z}}{p^n\mathbb{Z}}$ if $0 \leq n < \infty$ and $\widehat{\mathbb{Z}}(p^\infty) \stackrel{\text{def}}{=} \widehat{\mathbb{Z}}_p$ for $p \in \mathbb{P}$. With the conventions $p^\infty \widehat{\mathbb{Z}}_p \stackrel{\text{def}}{=} 0$ and $p^\infty \widehat{\mathbb{Z}}_q \stackrel{\text{def}}{=} \widehat{\mathbb{Z}}_q$ for $p \neq q$, we see that $p^n \widehat{\mathbb{Z}} = p^n \widehat{\mathbb{Z}}_p \times \prod_{p \neq q \in \mathbb{P}} \widehat{\mathbb{Z}}_q$ and $\frac{\widehat{\mathbb{Z}}}{p^n \widehat{\mathbb{Z}}} \cong \widehat{\mathbb{Z}}(p^n)$ for $p \in \mathbb{P}$ and $0 \leq n \in \mathbb{Z} \cup \{\infty\}$.

**Proposition 1.** *A nonzero finitely generated profinite abelian group is isomorphic to*

$$\Delta = \prod_{j=1}^{m} \prod_{p \in \mathbb{P}} \widehat{\mathbb{Z}}(p^{r_p(j)}) \text{ for some } 0 \leq r_p(j) \in \mathbb{Z} \cup \{\infty\},$$

$p \in \mathbb{P}, 1 \leq j \leq m$, *where* $r_p(j) \geq r_p(k) \Leftrightarrow j \leq k$, *and* $r_q(m) > 0$ *for some* $q \in \mathbb{P}$.

**Proof.** By Theorems 4.3.5 and 4.3.6 in [8], a finitely generated profinite abelian group is isomorphic to (*can be represented as*) $\prod_{j=1}^{m} \prod_{p \in \mathbb{P}} \widehat{\mathbb{Z}}(p^{r_p(j)})$ for some $0 \leq r_p(j) \in \mathbb{Z} \cup \{\infty\}, p \in \mathbb{P}, 1 \leq j \leq m$. The representation is indexed by $\{1, \ldots, m\} \times \mathbb{P}$. With regard to uniqueness up to isomorphism, there is no significance to the order of the factors $\widehat{\mathbb{Z}}(p^{r_p(j)})$ appearing. As long as the exact same aggregate list of $r_p(j)$ appears in such a representation, the associated groups are isomorphic.

For each $p \in \mathbb{P}$, we rearrange the $m$ exponents $r_p(1), \ldots, r_p(m)$ into descending order and relabel the ordered exponents $s_p(1), \ldots, s_p(m)$: $\{r_p(1), \ldots, r_p(m)\} = \{s_p(1), \ldots, s_p(m)\}$ where $s_p(1) \geq s_p(2) \geq \cdots \geq s_p(m)$. If, after applying this ordering for each $p \in \mathbb{P}$, we get $r_p(m) = 0$ for all $p \in \mathbb{P}$, then we remove all $\widehat{\mathbb{Z}}(p^{r_p(m)})$ for $p \in \mathbb{P}$, and reduce the value of $m$ accordingly. We repeat this weaning process right-to-left, so it terminates in a finite number of steps because $1 \leq m \in \mathbb{Z}$. In this way, we see that, without loss of generality, $m$ is minimal for a representation with the given characteristics. □

Define the *standard representation* of a nonzero finitely generated profinite abelian group to be the unique $\Delta = \prod_{j=1}^{m} \prod_{p \in \mathbb{P}} \widehat{\mathbb{Z}}(p^{r_p(j)})$ in Proposition 1 to which it is isomorphic. Set $\Delta_j \stackrel{\text{def}}{=} \prod_{p \in \mathbb{P}} \widehat{\mathbb{Z}}(p^{r_p(j)})$, $1 \leqslant j \leqslant m$, and $\Delta_p \stackrel{\text{def}}{=} \prod_{j=1}^{m} \widehat{\mathbb{Z}}(p^{r_p(j)})$, $p \in \mathbb{P}$.

Let $D$ be a finitely generated profinite abelian group. For $D \neq 0$ with standard representation $\Delta$ as in Proposition 1, define the *non-Archimedean width* of $D$ to be $\text{width}_{\text{nA}} D \stackrel{\text{def}}{=} m$ and set $\text{width}_{\text{nA}}\{0\} \stackrel{\text{def}}{=} 0$. Define the *non-Archimedean dimension* of $D$ to be $\dim_{\text{nA}} D \stackrel{\text{def}}{=} |\{j \in \{1, \ldots, m\} : \Delta_j \text{ is infinite}\}|$.

**Corollary 1.** *Non-Archimedean dimension of finitely generated profinite abelian groups is well-defined.*

**Proof.** Isomorphic finitely generated profinite groups have the same standard representation. □

A *supernatural vector* is any $\vec{\mathbf{n}} = (\mathbf{n}_1, \ldots, \mathbf{n}_m) \in \mathbb{S}^m$, $0 \leqslant m \in \mathbb{Z}$. Set $\mathbf{1} \stackrel{\text{def}}{=} \prod_{p \in \mathbb{P}} p^0 \in \mathbb{S}$ and $\vec{\mathbf{1}} \stackrel{\text{def}}{=} (\mathbf{1}, \mathbf{1}, \ldots, \mathbf{1}) \in \mathbb{S}^m$. Fix a finitely generated profinite abelian group $\Delta = \prod_{j=1}^{m} \prod_{p \in \mathbb{P}} \widehat{\mathbb{Z}}(p^{r_p(j)})$ as in Proposition 1. We write $\widehat{\mathbb{Z}}(\mathbf{n}) \stackrel{\text{def}}{=} \prod_{p \in \mathbb{P}} \widehat{\mathbb{Z}}(p^{\mathbf{n}_p})$ for $\mathbf{n} \in \mathbb{S}$ and $\widehat{\mathbb{Z}}(\vec{\mathbf{n}}) \stackrel{\text{def}}{=} \prod_{j=1}^{m} \prod_{p \in \mathbb{P}} \widehat{\mathbb{Z}}(p^{\mathbf{n}_{jp}})$ for $\vec{\mathbf{n}} \in \mathbb{S}^m$; note that supernatural vectors $\vec{\mathbf{n}}, \vec{\mathbf{n}}' \in \mathbb{S}^m$ associated to standard representations of isogenous finitely generated profinite groups agree in each coordinate except for a finite number of primes for which the exponents are finite. We introduce the notation $\mathbf{n}\widehat{\mathbb{Z}} \stackrel{\text{def}}{=} \prod_{p \in \mathbb{P}} p^{\mathbf{n}_p}\widehat{\mathbb{Z}}_p$ for $\mathbf{n} \in \mathbb{S}$ and $\vec{\mathbf{n}}\widehat{\mathbb{Z}}^m \stackrel{\text{def}}{=} \prod_{j=1}^{m} \prod_{p \in \mathbb{P}} p^{\mathbf{n}_{jp}}\widehat{\mathbb{Z}}_p$ for $\vec{\mathbf{n}} \in \mathbb{S}^m$.

**Corollary 2.** *A nonzero finitely generated profinite abelian group $\Delta$ is isomorphic to $\widehat{\mathbb{Z}}(\vec{\mathbf{n}})$ where $\vec{\mathbf{n}} \in \mathbb{S}^m$, $m = \text{width}_{nA} \Delta$, and the sequence $\vec{\mathbf{n}}\widehat{\mathbb{Z}}^m \rightarrowtail \widehat{\mathbb{Z}}^m \twoheadrightarrow \widehat{\mathbb{Z}}(\vec{\mathbf{n}})$ is exact.*

**Proof.** Proposition 1 gives that a finitely generated profinite abelian group is isomorphic to $\widehat{\mathbb{Z}}(\vec{\mathbf{n}})$ for some $\vec{\mathbf{n}} \in \mathbb{S}^m$. For each $p \in \mathbb{P}$, $\vec{\mathbf{n}}\widehat{\mathbb{Z}}^m$ has $p$-Sylow subgroup $\prod_{j=1}^{m} p^{\mathbf{n}_{jp}}\widehat{\mathbb{Z}}_p$. Because $\frac{\widehat{\mathbb{Z}}}{p^r\widehat{\mathbb{Z}}} \cong \widehat{\mathbb{Z}}(p^r)$ for $0 \leqslant r \in \mathbb{Z} \cup \{\infty\}$, we get $\frac{\widehat{\mathbb{Z}}^m}{\vec{\mathbf{n}}\widehat{\mathbb{Z}}^m} = \frac{\widehat{\mathbb{Z}}^m}{\prod_{j=1}^{m} \prod_{p \in \mathbb{P}} p^{\mathbf{n}_{jp}}\widehat{\mathbb{Z}}_p} \cong \prod_{j=1}^{m} \frac{\widehat{\mathbb{Z}}}{\prod_{p \in \mathbb{P}} p^{\mathbf{n}_{jp}}\widehat{\mathbb{Z}}_p} \cong \prod_{j=1}^{m} \frac{\prod_{p \in \mathbb{P}} \widehat{\mathbb{Z}}_p}{\prod_{p \in \mathbb{P}} p^{\mathbf{n}_{jp}}\widehat{\mathbb{Z}}_p} \cong \prod_{j=1}^{m} \prod_{p \in \mathbb{P}} \frac{\widehat{\mathbb{Z}}_p}{p^{\mathbf{n}_{jp}}\widehat{\mathbb{Z}}_p} = \widehat{\mathbb{Z}}(\vec{\mathbf{n}})$. □

## 3. Structure of Protori

For a torus-free protorus $G$ with profinite subgroup $\Delta$ inducing a torus quotient, we have by (Corollary 8.47, [2]) an accompanying injective morphism $\exp_G \colon \mathfrak{L}(G) \to G$ given by $\exp_G(r) = r(1)$. Set

- $\mathbb{Z}_\Delta \stackrel{\text{def}}{=} \Delta \cap \exp_G \mathfrak{L}(G)$,
- $\Gamma_\Delta \stackrel{\text{def}}{=} \{(\alpha, -\exp_G^{-1}\alpha) \colon \alpha \in \mathbb{Z}_\Delta\}$,
- $\pi_\Delta \colon \Delta \times \mathfrak{L}(G) \to \Delta$, the projection map onto $\Delta$,
- $\pi_\mathbb{R} \colon \Delta \times \mathfrak{L}(G) \to \mathfrak{L}(G)$, the projection map onto $\mathfrak{L}(G)$.

Then, $\pi_\Delta(\Gamma_\Delta) = \mathbb{Z}_\Delta$ and $\pi_\mathbb{R}(\Gamma_\Delta) = \exp_G^{-1} \mathbb{Z}_\Delta$ by the *Resolution Theorem for Compact Abelian Groups* (Theorem 8.20, [2]). *Note:* $\overline{K}$ for $K \subseteq G$ will always mean closure of $K$ in $G$ unless explicitly stated otherwise.

**Lemma 4.** *If $\Delta$ is a profinite subgroup of a torus-free protorus $G$ such that $G/\Delta \cong_t \mathbb{T}^{\dim G}$, then $\varphi_\Delta \colon \Delta \times \mathfrak{L}(G) \to G$, given by $\varphi_\Delta(\alpha, r) = \alpha + \exp_G r$, satisfies $\ker \varphi_\Delta \cong_t \mathbb{Z}^{\dim G}$.*

**Proof.** By (Theorem 8.20, [2]), $\ker \varphi_\Delta = \Gamma_\Delta$ and the projection $\pi_\mathbb{R} \colon \Delta \times \mathfrak{L}(G) \to \mathfrak{L}(G)$ restricts to a topological isomorphism $\pi_\mathbb{R}|_{\Gamma_\Delta} : \Gamma_\Delta \to \exp_G^{-1}\Delta = \exp_G^{-1}\mathbb{Z}_\Delta$, where $\exp_G$ is injective because $G$ is torus-free (Corollary 8.47, [2]). In addition, $\mathfrak{L}(G) \cong_t \mathbb{R}^{\dim G}$ and $\Gamma_\Delta$ is discrete by (Theorem 8.22 (6) $\Rightarrow$ (7), [2]). Thus, $\Gamma_\Delta \cong_t \exp_G^{-1}\mathbb{Z}_\Delta \cong_t \mathbb{Z}^k$ for some $0 \leq k \leq \dim G$ (Theorem A1.12, [2]). However, $[\Delta \times \mathfrak{L}(G)]/\Gamma_\Delta \cong_t G$ is compact, so $k = \dim G$. Thus, $\ker \varphi_\Delta = \Gamma_\Delta \cong_t \mathbb{Z}^{\dim G}$ as discrete groups. □

The next lemma identifies a simultaneously set-theoretic, topological, and algebraic property unique to profinite subgroups in a protorus which induce tori quotients.

**Lemma 5.** *If $\Delta$ is a profinite subgroup of a torus-free protorus $G$ such that $G/\Delta \cong_t \mathbb{T}^{\dim G}$, then $\overline{\mathbb{Z}_\Delta} = \Delta$ and $\mathbb{Z}_\Delta$ is closed in the subspace $\exp_G \mathfrak{L}(G)$.*

**Proof.** By (Theorem 8.20, [2]), a profinite subgroup $\Delta$ such that $G/\Delta \cong_t \mathbb{T}^{\dim G}$ always exists and for such a $\Delta$ we have $G \cong_t G_\Delta \stackrel{\text{def}}{=} \frac{\Delta \times \mathfrak{L}(G)}{\Gamma_\Delta}$ where $\Gamma_\Delta = \{(\exp_G r, -r) : r \in \mathfrak{L}(G), \exp_G r \in \Delta\}$ is a free abelian group and rank $\Gamma_\Delta = \dim G = \text{rank}\,[\Delta \cap \exp_G \mathfrak{L}(G)]$ by Lemma 4 and the fact that $\exp_G$ is injective when $G$ is torus-free (Corollary 8.47, [2]). We have $\pi_\Delta(\Gamma_\Delta) = \mathbb{Z}_\Delta \subseteq \Delta' \stackrel{\text{def}}{=} \overline{\mathbb{Z}_\Delta}$, so $\Gamma_\Delta$ is a subgroup of $\Delta' \times \mathfrak{L}(G)$. Because $\Delta \rightarrowtail \frac{\Delta \times \{0\} + \Gamma_\Delta}{\Gamma_\Delta} \subset G_\Delta$ is a topological isomorphism onto its image, $\Delta' \rightarrowtail \frac{\Delta' \times \{0\} + \Gamma_\Delta}{\Gamma_\Delta} \subset G_\Delta$ is as well. Since $\Gamma_\Delta$ is discrete in $\Delta \times \mathfrak{L}(G)$ (Theorem 8.20, [2]), it is discrete in $\Delta' \times \mathfrak{L}(G)$, so $G_{\Delta'} \stackrel{\text{def}}{=} [\Delta' \times \mathfrak{L}(G)]/\Gamma_\Delta$ is a Hausdorff subgroup of $G_\Delta$. However, $(\Delta \setminus \Delta') \times \mathfrak{L}(G)$ is open in $\Delta \times \mathfrak{L}(G)$ and the quotient map $q_\Delta \colon \Delta \times \mathfrak{L}(G) \to G_\Delta$ is an open map, so $q_\Delta[(\Delta \setminus \Delta') \times \mathfrak{L}(G)] = [(\Delta \setminus \Delta') \times \mathfrak{L}(G) + \Gamma_\Delta]/\Gamma_\Delta = G_\Delta \setminus G_{\Delta'}$, is open in $G_\Delta$. It follows that $G_{\Delta'}$ is a compact abelian subgroup of $G_\Delta$ and $\frac{G_\Delta}{G_{\Delta'}} = \frac{[\Delta \times \mathfrak{L}(G)]/\Gamma_\Delta}{[\Delta' \times \mathfrak{L}(G)]/\Gamma_\Delta} \cong_t \frac{\Delta \times \mathfrak{L}(G)}{\Delta' \times \mathfrak{L}(G)} \cong_t \frac{\Delta}{\Delta'}$ by (Theorem 5.35, [5]). Thus, there is an exact sequence $G_{\Delta'} \rightarrowtail G_\Delta \twoheadrightarrow \frac{\Delta}{\Delta'}$. Now, $\dim \Delta = 0 \Rightarrow \dim(\Delta/\Delta') = 0$ and we know $\Delta/\Delta'$ is compact Hausdorff, so $\Delta/\Delta'$ is totally disconnected (Corollary 7.72, [4]). Thus, $(\Delta/\Delta')^\vee$ is torsion (Corollary 8.5, [2]). By Pontryagin duality, $(\Delta/\Delta')^\vee$ embeds in the torsion-free group $G_\Delta^\vee$, so $(\Delta/\Delta')^\vee = 0$ and $\Delta = \Delta' = \overline{\mathbb{Z}_\Delta}$.

Lastly, if $x$ lies in the closure of $\mathbb{Z}_\Delta$ in $\exp_G \mathfrak{L}(G)$ under the (metric) subspace topology, then $x \in \exp_G \mathfrak{L}(G)$ and $x$ is the limit of a sequence of elements of $\mathbb{Z}_\Delta$. However, $\Delta$ is closed, so $x \in \Delta \cap \exp_G \mathfrak{L}(G) = \mathbb{Z}_\Delta$. This proves that $\mathbb{Z}_\Delta$ is closed in the subspace $\exp_G \mathfrak{L}(G)$. □

Define $\mathrm{L}(G) = \{\Delta \subset G \colon \Delta \text{ a profinite subgroup such that } G/\Delta \text{ is a torus}\}$ for a protorus $G$. In Proposition 2, we show that $\mathrm{L}(G)$ is a lattice with $+$ as join and $\cap$ as meet; in particular, $\mathrm{L}(G)$ is directed upward and downward. We then prove a number of useful closure properties for $\mathrm{L}(G)$, which are applied for the remainder. However, first a remark regarding some facts which we use freely going forward without further mention.

**Remark 1.** (i) *By (Theorem 8.46.(iii), [2]), a path-connected protorus is a torus. Thus, if a protorus $G$ has a profinite subgroup $D$ and $G/D$ is the continuous image of a torus, then automatically $D \in \mathrm{L}(G)$. (ii) A profinite subgroup $D$ of a torus $H$ is finite: $H/D$ is a path-connected protorus, whence a torus of the same dimension as $H$ by Lemma 1. By duality, $(H/D)^\vee$ and $H^\vee$ are free abelian with $\mathrm{rk}(H/D)^\vee = \mathrm{rk}\, H^\vee$ and $(H/D)^\vee \rightarrowtail H^\vee \twoheadrightarrow D^\vee$ exact, so $D^\vee$, whence $D$, is finite.*

**Proposition 2.** *For a torus-free protorus $G$, $\mathrm{L}(G)$ is a countable lattice under $\cap$ for meet and $+$ for join. $\mathrm{L}(G)$ is closed under:*

1. *preimages via $\mu_n$, $0 \neq n \in \mathbb{Z}$.*
2. *finite extensions.*
3. *scalar multiplication by nonzero integers.*
4. *join ($+$). and*
5. *meet ($\cap$).*

*Given any $\Delta, \Delta' \in L(G)$, there exists $0 < k \in \mathbb{Z}$ such that $k\Delta \subseteq \Delta'$. If $\Delta' \subseteq \Delta$, then $[\Delta \colon \Delta'] < \infty$.*

**Proof.** Each $\Delta \in L(G)$ corresponds via Pontryagin duality to a unique-up-to-isomorphism torsion abelian quotient of $X = G^\vee$ by a free abelian subgroup $Z_\Delta$ with $\mathrm{rk} Z_\Delta = \mathrm{rk} X$. Because $X$ is countable and there are countably many finite subsets of a countable set (corresponding to bases of $Z_\Delta$s, counting one basis per $Z_\Delta$), it follows that $L(G)$ is countable.

1. $\mu_n \colon G \to G$ has finite kernel by Lemma 2 so its restriction $\mu_n^{-1}\Delta \to \Delta$ has finite kernel for $\Delta \in L(G)$. Since $\ker \mu_n$ and $\Delta \in L(G)$ are zero-dimensional compact abelian groups, it follows from Lemma 1 that the compact Hausdorff subgroup $\mu_n^{-1}\Delta$ is zero-dimensional, whence profinite. Because the natural map $G/\Delta \to G/\mu_n^{-1}\Delta$ is surjective and $G/\Delta$ is a torus, it follows that $G/\mu_n^{-1}\Delta$ is path-connected, whence $G/\mu_n^{-1}\Delta$ is a torus (Theorem 8.46.(iii), [2]) and $\mu_n^{-1}\Delta \in L(G)$.

2. If $\Delta \in L(G)$ has index $1 \leq m \in \mathbb{Z}$ in a subgroup $D$ of $G$, then $D$ is the sum of finitely many copies of $\Delta$, so $D$ is compact. Thus, $\Delta \subseteq D \subseteq \mu_m^{-1}\Delta \in L(G)$ by 1, so $D$ is profinite. The natural morphism $G/\Delta \to G/D$ is surjective and $G/\Delta$ a torus, so $D \in L(G)$.

3. $\mu_j|_\Delta \colon \Delta \to j\Delta$ is surjective with finite kernel by Lemma 2, so $j\Delta$ is profinite. $G$ is divisible so $\mu_j \colon G \to G$ is surjective, thus inducing a surjective morphism $\frac{G}{\Delta} \to \frac{G}{j\Delta}$. It follows that $j\Delta \in L(G)$.

4. Addition defines a surjective morphism $\Delta \times \Delta' \twoheadrightarrow \Delta + \Delta'$. By Lemma 1, it follows that $\Delta \times \Delta'$, whence $\Delta + \Delta'$, is profinite. Because the natural map $\frac{G}{\Delta} \to \frac{G}{\Delta+\Delta'}$ is surjective, $\Delta + \Delta' \in L(G)$.

5. The kernel of $\frac{G}{\Delta} \twoheadrightarrow \frac{G}{\Delta+\Delta'}$ is $\frac{\Delta+\Delta'}{\Delta}$, a zero-dimensional subgroup of $\frac{G}{\Delta}$ by Lemma 1. As a zero-dimensional subgroup of a torus, $\frac{\Delta+\Delta'}{\Delta} \cong \frac{\Delta'}{\Delta \cap \Delta'}$ is finite, so there is a nonzero integer $l$ such that $l\Delta' \subseteq \Delta$. Lemma 1 gives that $\Delta \cap \Delta'$ is zero-dimensional, whence profinite. We know by 3 that $l\Delta' \in L(G)$, thus the natural map $\frac{G}{l\Delta'} \to \frac{G}{\Delta \cap \Delta'}$ is a surjective morphism, whence $\Delta \cap \Delta' \in L(G)$.

It follows from 4 and 5 that $L(G)$ is a lattice. It remains to show that, if $\Delta' \subseteq \Delta$, then $[\Delta \colon \Delta'] < \infty$. Arguing as in 5, there exists $0 < k \in \mathbb{Z}$ such that $k\Delta \subseteq \Delta'$. $\Delta$ is a finitely generated profinite abelian group, thus $[\Delta \colon \Delta'] \leq [\Delta \colon k\Delta] < \infty$. □

**Corollary 3.** *Elements of $L(G)$ are mutually isogenous in a torus-free protorus $G$.*

**Proof.** Suppose that $\Delta_1, \Delta_2 \in L(G)$. We proved in Proposition 2 that there exist nonzero integers $k$ and $l$ such that $k\Delta_1 \subseteq \Delta_2$, $l\Delta_2 \subseteq \Delta_1$, $[\Delta_2 \colon k\Delta_1] < \infty$, and $[\Delta_1 \colon l\Delta_2] < \infty$. The multiplication-by-$k$ and multiplication-by-$l$ morphisms thus exhibit an isogeny between $\Delta_1$ and $\Delta_2$. □

**Lemma 6.** *Non-Archimedean dimension of finitely generated profinite groups is invariant under isogeny.*

**Proof.** If two such groups, say $D$ and $D'$ are isogenous, then so are their standard representations, say $\Delta_D$ and $\Delta_{D'}$, as in Proposition 1. Multiplying both groups by the same sufficiently large integer, say $N$, produces isogenous groups $ND$ and $ND'$ with $\dim_{\mathrm{nA}}(ND) = \mathrm{width}_{\mathrm{nA}}(ND)$ and $\dim_{\mathrm{nA}}(ND') = \mathrm{width}_{\mathrm{nA}}(ND')$. By Corollary 2, $ND$ and $ND'$ have standard representations, say $\widehat{\mathbb{Z}}(\vec{\mathbf{n}})$ and $\widehat{\mathbb{Z}}(\vec{\mathbf{n}}')$ for some supernatural vectors $\vec{\mathbf{n}}$ and $\vec{\mathbf{n}}'$. If $\widehat{\mathbb{Z}}(\vec{\mathbf{n}})$ and $\widehat{\mathbb{Z}}(\vec{\mathbf{n}}')$ have distinct non-Archimedean dimensions, then one of the two has an *extra* coordinate $k$ for which there are factors $\widehat{\mathbb{Z}}(p^{n_{kp}})$ for infinitely many primes $p$ and/or factors isomorphic to one or more copies of some $\widehat{\mathbb{Z}}_q$ for some prime $q$; this is impossible if $\widehat{\mathbb{Z}}(\vec{\mathbf{n}})$ and $\widehat{\mathbb{Z}}(\vec{\mathbf{n}}')$ are isogenous because supernatural vectors associated to standard representations of isogenous groups can differ for any given coordinate at only a finite number of primes and only at those primes having finite exponents. Thus, the definition of non-Archimedean dimension and its preservation under multiplication by $N$ give that $\dim_{\mathrm{nA}} D = \dim_{\mathrm{nA}}(ND) = \dim_{\mathrm{nA}} \widehat{\mathbb{Z}}(\vec{\mathbf{n}}) = \dim_{\mathrm{nA}} \widehat{\mathbb{Z}}(\vec{\mathbf{n}}') = \dim_{\mathrm{nA}}(ND') = \dim_{\mathrm{nA}} D'$. □

Define the ***non-Archimedean dimension*** of a protorus $G$ to be $\dim_{\mathrm{nA}} G \stackrel{\mathrm{def}}{=} \dim_{\mathrm{nA}} \Delta$ for a profinite subgroup $\Delta$ of $G$ for which $G/\Delta$ is a torus.

**Corollary 4.** *Non-Archimedean dimension of protori is well-defined.*

**Proof.** Elements of L(*G*) for a protorus *G* are isogenous by Corollary 3, so the result follows by Lemma 6. □

A protorus *G* with profinite subgroup $\Delta$ satisfying $G/\Delta \cong_t \mathbb{T}^{\dim G}$ will always satisfy $\dim_{nA} G = \dim_{nA} \Delta \leq \text{width}_{nA} \Delta \leq \dim G$; and for *G* torus-free, *always* $\text{rk}\,\mathbb{Z}_\Delta = \dim G$.

A protorus *G* is ***factorable*** if there are non-trivial protori $G_1$ and $G_2$ such that $G \cong_t G_1 \times G_2$, and *G* is ***completely factorable*** if $G \cong_t \prod_{i=1}^n G_i$ where $\dim G_i = 1, 1 \leq i \leq n$. A result by Mader and Schultz [10] has the surprising implication that the classification of protori up to topological near-isomorphism reduces to that of protori with no one-dimensional factors.

**Proposition 3.** *If D is a finitely generated profinite abelian group, then there is a completely factorable protorus G containing a closed subgroup $\Delta \cong D$ such that $G/\Delta$ is a torus.*

**Proof.** First, note that the finite cyclic group $\mathbb{Z}(r)$, $0 < r \in \mathbb{Z}$, is isomorphic to the closed subgroup $\frac{(1/r)\mathbb{Z}}{\mathbb{Z}}$ of the torus $\frac{\mathbb{R}}{\mathbb{Z}}$, so it follows that $\frac{(1/r)\mathbb{Z}}{\mathbb{Z}}$ is a profinite subgroup of $\frac{\mathbb{R}}{\mathbb{Z}}$ inducing a torus quotient. Next, by Proposition 1 there is no loss of generality in assuming $D = \widehat{\mathbb{Z}}(\vec{n})$ for some $\vec{n} \in \mathbb{S}^m$ where $\Delta_j(\vec{n}) \neq 0$ for $1 \leq j \leq m$. If $\Delta_j(\vec{n})$ is finite then it must be isomorphic to $\mathbb{Z}(r_j)$ for some $0 < r_j \in \mathbb{Z}$; in this case, set $G_j \stackrel{\text{def}}{=} \frac{\mathbb{R}}{\mathbb{Z}}$ and $E_j \stackrel{\text{def}}{=} \frac{(1/r_j)\mathbb{Z}}{\mathbb{Z}}$. If $\Delta_j(\vec{n})$ is not finite, then $G_j \stackrel{\text{def}}{=} [\Delta_j(\vec{n}) \times \mathbb{R}]/\mathbb{Z}(\mathbf{1}, 1)$ is a solenoid (1-dimensional protorus) containing a closed subgroup $E_j \cong_t \Delta_j(\vec{n})$ satisfying $G_j/E_j \cong_t \mathbb{T}$ (Theorem 10.13, [5]). It follows that $G \stackrel{\text{def}}{=} G_1 \times \cdots \times G_m$ is a protorus containing the closed subgroup $\Delta \stackrel{\text{def}}{=} E_1 \times \cdots \times E_m \cong_t D$ and satisfying $G/\Delta \cong_t \mathbb{T}^m$. □

**Corollary 5.** *If D is a finitely generated profinite abelian group with $\text{width}_{nA}\, D = \dim_{nA}\, D$, then there is a completely factorable torus-free protorus G containing a closed subgroup $\Delta \cong D$ such that $G/\Delta$ is a torus.*

**Proof.** In this case, no $\Delta_j(\mathbf{n})$ factor is finite cyclic in the proof of Proposition 3. □

An arbitrary topological group *K* is ***compactly ruled*** if it is a directed union of compact open subgroups (Definition 1.4, [3]). An abelian compactly ruled group satisfies $\overline{\mathbb{Z}g}$ is compact for each $g \in K$ (Corollary 1.12, [3]). A locally compact abelian group *K* is ***periodic*** if it is totally disconnected and $\overline{\mathbb{Z}g}$ is compact for each $g \in K$ (Definition 1.13, [3]). Abbreviate *totally disconnected locally compact* **tdlc** and locally compact abelian **LCA**. Every element *g* in a periodic LCA group *K* satisfies $\overline{\mathbb{Z}g} \cong \prod_{p \in \mathbb{P}} \widehat{\mathbb{Z}}(p^{s_p})$ for some $s \in \mathbb{S}$ and there corresponds a $\widehat{\mathbb{Z}}$-module structure on *K* with continuous scalar multiplication $\widehat{\mathbb{Z}} \times K \to K$ compatible with the unique *p*-Sylow decomposition of *K* (pp. 48–49, [3]): $z \cdot h \stackrel{\text{def}}{=} (z \cdot h)_{p \in \mathbb{P}} = (z)_{p \in \mathbb{P}} \cdot (h)_{p \in \mathbb{P}}$ for $z \in \widehat{\mathbb{Z}}$ and $h \in \overline{\mathbb{Z}g}$ (Proposition 4.21, [3]). A periodic LCA group $\widehat{\Delta}$ is a ***topological divisible hull*** of a profinite group $\Delta$ if $\widehat{\Delta}$ is an algebraic divisible hull of $\Delta$ with the unique group topology for which an isomorphic copy of $\Delta$ is an open subgroup (Theorem 3.42, [3]).

For a torus-free protorus *G* and a sublattice $L \subseteq L(G)$ satisfying $\bigcap L = 0$, set

- $\widehat{\Delta}_L \stackrel{\text{def}}{=} \sum_{\Delta \in L} \Delta$.
- $V_L \stackrel{\text{def}}{=} \sum_{\Delta \in L} \mathbb{Z}_\Delta = \bigcup_{\Delta \in L} \mathbb{Z}_\Delta$.
- $\Gamma_L \stackrel{\text{def}}{=} \{(\alpha, -\exp_G^{-1} \alpha) : \alpha \in V_L\}$.

**Proposition 4.** *If G is a torus-free protorus and L is a sublattice of L(G) satisfying $\bigcap L = 0$, then the group $\widehat{\Delta}_L$ topologized by taking L to be an open neighborhood basis at 0 is a periodic $\widehat{\mathbb{Z}}$-module with $\widehat{\mathbb{Z}}g = \overline{\mathbb{Z}g}$ for each $g \in \widehat{\Delta}_L$, and $\widehat{\Delta}_L \times \mathfrak{L}(G)$ is locally compact with closed subgroup $\Gamma_L$. The identity map to the subspace topology on $\widehat{\Delta}_L$ is continuous.*

**Proof.** L($G$) is closed under finite sums and directed upward by Proposition 2, so $\widehat{\Delta}_L = \bigcup_{\Delta \in L} \Delta$. The elements of $L$ are finitely generated profinite abelian groups with $[D_i : D_1 \cap D_2] < \infty$ for $i = 1, 2$ for any $D_1, D_2 \in L$, and thus are compact in the topology defined on $\widehat{\Delta}_L$ by declaring $L$ to be a neighborhood basis of 0. This Hausdorff topology ($\bigcap L = 0$) is therefore compactly ruled, whence locally compact (Proposition 1.3, [3]). It follows that $\widehat{\Delta}_L \times \mathfrak{L}(G)$ is locally compact. By Proposition 2, $L$ is a countable neighborhood basis at 0 of compact open subgroups, and so $\widehat{\Delta}_L$ is metrizable (Theorem 8.3, [5]) and totally disconnected (Theorem 1.34, [2]). Thus, $\widehat{\Delta}_L$ is periodic, whence $\widehat{\mathbb{Z}}g = \overline{\mathbb{Z}g}$ for each $g \in \Delta \in L$ and $\widehat{\Delta}_L$ is a $\widehat{\mathbb{Z}}$-module with continuous bilinear scalar multiplication (Proposition 4.22, [3]). To see that $\Gamma_L$ is closed, suppose that $\{(\alpha_j, r_j)\}_{j \geq 1} \subset \Gamma_L$ converges to some $(\beta, s) \in \widehat{\Delta}_L \times \mathfrak{L}(G)$. Then, $\beta$ and $\{\alpha_j\}_{j \geq N}$ lie in some $\Delta_0 \in L$ for some $1 \leq N \in \mathbb{Z}$. We can assume without loss of generality that $N = 1$. Since $r_j \to s \in \mathfrak{L}(G)$, continuity of $\exp_G$ implies $-\alpha_j = \exp_G(r_j) \to \exp_G(s) = -\beta \in \Delta_0 \cap \exp_G \mathfrak{L}(G) = \mathbb{Z}_{\Delta_0} \subset V_L$, whence $(\beta, s) \in \Gamma_L$. Next, by (Theorem 8.22, [2]), a basic open neighborhood of $G$ can be taken to have the form $\Delta + \exp_G B$ for $\Delta \in L$ and $B$ a Euclidean open ball about 0 in $\mathfrak{L}(G)$; then $\Delta$ is an open subset of $\widehat{\Delta}_L$ under its locally compact topology, and $\Delta$ is contained in the open subset $(\Delta + \exp_G B) \cap \widehat{\Delta}_L$ of $\widehat{\Delta}_L$ under its subspace topology, so the identity map is continuous. □

The group of path components of a protorus $G$ is $\pi_0(G) \stackrel{\text{def}}{=} G / \exp_G \mathfrak{L}(G)$. The next proposition describes properties of substructures of a protorus, to be fine-tuned in Theorem 1.

**Proposition 5.** *If $G$ is a torus-free protorus with no factors $\cong_t \mathbb{Q}^\vee$, then*

1. tor $G = \text{tor } \widehat{\Delta}_{L(G)}$, $\widehat{\Delta}_{L(G)}$ *is divisible, and $V_{L(G)}$ is dense in $G$.*
2. $\widehat{\Delta}_{L(G)}/\Delta \cong (\mathbb{Q}/\mathbb{Z})^{\dim G}$ *for each $\Delta \in L(G)$ and $G/\widehat{\Delta}_{L(G)} \cong (\mathbb{R}/\mathbb{Q})^{\dim G}$.*
3. $\widehat{\Delta}_{L(G)} \cong_t \varprojlim_{\Delta \in L(G)} (\widehat{\Delta}_{L(G)}/\Delta)$ *where $\widehat{\Delta}_{L(G)}/\Delta$ is a discrete group for each $\Delta \in L(G)$.*
4. $\exp_G \mathfrak{L}(G)$, *the path component 0 in $G$, is divisible and torsion-free.*
5. $\pi_0(G) \cong \text{Ext}(G^\vee, \mathbb{Z}) \cong \Delta/\mathbb{Z}_\Delta$ *for each $\Delta \in L(G)$.*

**Proof.** 1. If $x \in G$ and $nx \in \widehat{\Delta}_{L(G)}$ for some $0 \neq n \in \mathbb{Z}$, then $nx \in \Delta$ for some $\Delta \in L(G)$, so $x \in \mu_n^{-1}\Delta \subset \widehat{\Delta}_{L(G)}$ by Proposition 2. Thus, $G/\widehat{\Delta}_{L(G)}$ is torsion-free whence tor $G = \text{tor } \widehat{\Delta}_{L(G)}$. If $y \in \widehat{\Delta}_{L(G)}$ and $0 \neq k \in \mathbb{Z}$, then because $G$ is divisible there is a $z \in G$ such that $kz = y$; as we have shown this implies $z \in \widehat{\Delta}_{L(G)}$; this proves that $\widehat{\Delta}_{L(G)}$ is divisible. Because tor $G$ is dense in $G$ (Corollary 8.9, [2]), for $V_{L(G)}$ to be dense in $G$ it suffices to show that $V_{L(G)}$ is dense in $\widehat{\Delta}_{L(G)}$: $\overline{V_{L(G)}} = \overline{\sum_{\Delta \in L} \mathbb{Z}_\Delta} \supseteq \sum_{\Delta \in L} \overline{\mathbb{Z}_\Delta} = \sum_{\Delta \in L} \Delta = \widehat{\Delta}_{L(G)}$.

2. Let $\Delta \in L(G)$. By 1, $G/\widehat{\Delta}_{L(G)} \cong (G/\Delta)/(\widehat{\Delta}_{L(G)}/\Delta)$ is torsion-free so $G/\Delta \cong \mathbb{T}^{\dim G}$, whence $\text{tor}(G/\Delta) \subset \widehat{\Delta}_{L(G)}/\Delta$. Conversely, if $x \in \widehat{\Delta}_{L(G)}$, then $x \in \Delta'$ for some $\Delta' \in L(G)$, whence $kx \in \Delta$ for some $0 \neq k \in \mathbb{Z}$ by Proposition 2; that is, $x + \Delta \in \text{tor}(\widehat{\Delta}_{L(G)}/\Delta)$. Thus, $\widehat{\Delta}_{L(G)}/\Delta = \text{tor}(G/\Delta) \cong \text{tor}(\mathbb{T}^{\dim G}) = (\mathbb{Q}/\mathbb{Z})^{\dim G}$ so $G/\widehat{\Delta}_{L(G)} \cong \mathbb{T}^{\dim G}/(\mathbb{Q}/\mathbb{Z})^{\dim G} = (\mathbb{R}/\mathbb{Z})^{\dim G}/(\mathbb{Q}/\mathbb{Z})^{\dim G} \cong (\mathbb{R}/\mathbb{Q})^{\dim G}$.

3. Under its locally compact topology, $\widehat{\Delta}_{L(G)}$ is periodic (Proposition 4), so the result follows by (pp. 48–49, [3]).

4. Theorem 8.30 [2], gives that $\exp_G \mathfrak{L}(G)$ is the path component of 0. $G$ is torus-free, so the corestriction $\exp_G \colon \mathfrak{L}(G) \to \exp_G \mathfrak{L}(G)$ is an isomorphism of abelian groups (Corollary 8.47, [2]). Thus, $\exp_G \mathfrak{L}(G)$ is divisible and torsion-free because $\mathfrak{L}(G) \cong_t \mathbb{R}^{\dim G}$ is (Propositions 7.25 and 7.36, [2]).

5. (Theorem 8.30, [2]) and (Corollary 8.33, [2]). □

We define an apparatus upon which our proof of Theorem 1 depends. The setting involves an $n$-dimensional torus-free protorus $G$. For such $G$, $\exp_G \colon \mathfrak{L}(G) \to G$ is injective (Corollary 8.47, [2]) and $\mathbb{R}^n \cong_t \mathfrak{L}(G)$ (Proposition 7.24, [2]). We argue relative to a fixed $\Delta^* \in L(G)$ and its dense rank-$n$ free abelian subgroup $\mathbb{Z}_{\Delta^*} = \Delta^* \cap \exp_G \mathfrak{L}(G)$ (Lemma 5). There exists an algebraic isomorphism $\theta_{\Delta^*} \colon \mathbb{R}^n \to \exp_G \mathfrak{L}(G)$ with $\mathbb{Z}_{\Delta^*} = \theta_{\Delta^*}(\mathbb{Z}^n) \subset \theta_{\Delta^*}(\mathbb{Q}^n) \subset \exp_G \mathfrak{L}(G) = \theta_{\Delta^*}(\mathbb{R}^n)$, where $\theta_{\Delta^*}(\mathbb{Q}^n) = \{w \in \exp_G \mathfrak{L}(G) \colon kw \in \Delta^* \text{ for some } 0 \neq k \in \mathbb{Z}\}$.

We fix a unit hemisphere $H$ and introduce a parameter $Y$ to represent an arbitrary rank-$n$ torsion-free abelian group with $H \subset \mathbb{Z}^n \subseteq Y \subseteq \mathbb{Q}^n$; note that $\theta_{\Delta^*}(h)/p \notin \mathbb{Z}_{\Delta^*}$ for any $h \in H$ and $p \in \mathbb{P}$. With $\Delta^* \in L(G)$, $\theta_{\Delta^*} \colon \mathbb{R}^n \to \exp_G \mathfrak{L}(G)$, and a unit hemisphere $H \subset \mathbb{Z}^n$ fixed, we itemize pertinent background information for Theorem 1:

(a) $X$ is a fixed value of $Y$ for which $\mathbb{Z}^n \rightarrowtail X \twoheadrightarrow X/\mathbb{Z}^n$ dualizes $\Delta^* \rightarrowtail G \twoheadrightarrow G/\Delta^*$.

(b) $\Delta_y \stackrel{\text{def}}{=} \Delta^* + \mathbb{Z}\theta_{\Delta^*}(y) \in L(G)$ and $\mathbb{Z}_{\Delta_y} = \mathbb{Z}_{\Delta^*} + \mathbb{Z}\theta_{\Delta^*}(y)$, $y \in Y$ (Proposition 2), $\widehat{\Delta}_Y \stackrel{\text{def}}{=} \sum_{y \in Y} \Delta_y$.

(c) $Y^* \stackrel{\text{def}}{=} \theta_{\Delta^*}(Y) = \sum_{y \in Y} \mathbb{Z}_{\Delta_y}$ and $L_Y \stackrel{\text{def}}{=} \{\sum_{y \in Y'} \Delta_y \colon Y' \subset Y, Y' \text{ finite}\} \cup \{\Delta \in L(G) \colon \Delta \subseteq \Delta^*\} \subseteq L(G)$.

(d) $s_z^Y \in \mathbb{S}$ where $\{s_z^Y(p)\}_{p \in \mathbb{P}}$ is the height sequence in $Y$ of $z \in H$.

(e) $s_{\text{fin}\infty}^Y(z) \stackrel{\text{def}}{=} \prod_{0 < s_z^Y(p) < \infty} p^\infty$, $z_{\text{fin}}^\infty \stackrel{\text{def}}{=} \{\frac{jz}{\ell} \colon j, \ell \in \mathbb{Z}, \ell \mid s_{\text{fin}\infty}^Y(z)\}$, and $Y_{\text{fin}}^\infty \stackrel{\text{def}}{=} \sum_{z \in H} z_{\text{fin}}^\infty$.

(f) $s_\infty^Y(z) \stackrel{\text{def}}{=} \prod_{p \mid s_z^Y} p^\infty$, $z_\infty \stackrel{\text{def}}{=} \{\frac{jz}{\ell} \colon j, \ell \in \mathbb{Z}, \ell \mid s_\infty^Y(z)\}$, and $Y_\infty \stackrel{\text{def}}{=} \sum_{z \in H} z_\infty$.

(g) $Y_{\text{fin}}^\infty$ and $Y_\infty$ are quotient-divisible: $Y_{\text{fin}}^\infty/\mathbb{Z}^n$ and $Y_\infty/\mathbb{Z}^n$ are divisible.

**Theorem 1.** (Structure Theorem for Protori) *A protorus is topologically isomorphic to $\mathbb{T}^r \times (\mathbb{Q}^\vee)^k \times G$ for some n-dimensional protorus $G$ with no factors $\cong_t \mathbb{T}$ or $\mathbb{Q}^\vee$ and nonnegative integers r,k,n where the three factors are uniquely determined up to topological isomorphism. In view of the definitions and notation above— with fixed $\Delta^*$, fixed unit hemisphere H, fixed $X \cong G^\vee$, and arbitrary Y with $H \subset \mathbb{Z}^n \subseteq X, Y \subseteq \mathbb{Q}^n$—the structure of G is as follows:*

1. $L_Y$ is a lattice, $V_{L_Y} = Y^* = \widehat{\Delta}_Y \cap \exp_G \mathfrak{L}(G)$, and $\Gamma_{L_Y} \cong Y$.
2. $\widehat{\Delta}_{L_Y} = \widehat{\Delta}_Y$ has a periodic LCA topology, and $Y^*$ is dense in $G \Leftrightarrow Y$ has no free summands.
3. $G \cong_t \frac{\widehat{\Delta}_Y \times \mathfrak{L}(G)}{Y}$; in particular, $G \cong_t \frac{\widehat{\Delta}_X \times \mathfrak{L}(G)}{X}$ where $X^* \cong X \cong G^\vee$ is dense in $G$.
4. $\widehat{\Delta}_{X_\infty}$ is a topological divisible hull of $\Delta^*$ and $\text{tor } G = \text{tor } \widehat{\Delta}_{L(G)} = \text{tor } \widehat{\Delta}_{\mathbb{Q}^n}$.
5. $\widehat{\Delta}_{X_\infty} \cong_t \prod_{p \in \mathbb{P}}^{\text{loc}}((\widehat{\Delta}_{X_\infty})_p, \Delta_p^*) \cong_t \varprojlim_{\Delta \in L_{X_\infty}}(\widehat{\Delta}_{X_\infty}/\Delta)$; $X^\vee \cong_t \varprojlim_{\Delta \in L_X}(G/\Delta)$; and $X \cong \varinjlim_{\Delta \in L_X}(G/\Delta)^\vee$.

**Proof.** Let $K$ be a protorus. By Corollary 3.8.3 in [9], $K^\vee = Z \oplus C$ where $Z$ is free abelian and $C$ is a subgroup with no free summands which is uniquely determined by $K^\vee$. By Corollary 4.2.5 in [9], $C = R \oplus D$ where $D$ is a torsion-free divisible group uniquely determined by $C$, and $R$ is unique up to isomorphism. The first statement of the theorem thus follows by duality. The remainder of the proof involves the $n$-dimensional protorus $G$ with no factors $\cong_t \mathbb{T}$ or $\mathbb{Q}^\vee$.

1. $Y^*$ is the directed union of free abelian groups $F$ with $\mathbb{Z}_{\Delta^*} \subseteq F \subseteq Y^*$, and $F/\mathbb{Z}_{\Delta^*}$ is finite for such $F$. Hence, $\Delta_F \stackrel{\text{def}}{=} \overline{F} \in L(G)$ by Proposition 2, whence $\Delta_F \cap \exp_G \mathfrak{L}(G) = \mathbb{Z}_{\Delta_F} = F$ (Lemma 5) and $\Delta_F \in L_Y$ for such $F$. Thus, $L_Y = \{\Delta_F \colon \mathbb{Z}_{\Delta^*} \subseteq F \subseteq Y, F \text{ free abelian}\} \cup \{\Delta \in L(G) \colon \Delta \subseteq \Delta^*\}$ is a lattice.

In addition, $\widehat{\Delta}_Y \cap \exp_G \mathfrak{L}(G) = (\sum_{y \in Y} \Delta_y) \cap \exp_G \mathfrak{L}(G) = \bigcup_{y \in Y}(\Delta_y \cap \exp_G \mathfrak{L}(G)) = \bigcup_{y \in Y} \mathbb{Z}_{\Delta_y} = \sum_{y \in Y} \mathbb{Z}_{\Delta_y} = Y^*$. By Lemma 4, $Y^* = V_{L_Y}$ and the closed subgroup $\Gamma_{L_Y} = \{(\alpha, -\exp_G^{-1}\alpha) \colon \alpha \in V_{L_Y}\}$ is equal to $\{(\alpha, -\exp_G^{-1}\alpha) \colon \alpha \in Y^*\} \cong Y^* \cong Y$.

2. By definition, $\widehat{\Delta}_{L_Y} = \sum_{\Delta \in L_Y} \Delta = \sum_{y \in Y} \Delta_y = \widehat{\Delta}_Y$, and $\bigcap_{\Delta \in L_Y} \Delta = \{0\}$, so $\widehat{\Delta}_Y$ has a periodic LCA topology by Proposition 4. The morphism $\exp_G$ is continuous and injective (Corollary 8.47, [2]) with $\exp_G \mathfrak{L}(G)$ dense in $G$. Thus, $Y^*$ is dense in $G \Leftrightarrow Y^*$ is dense in $\exp_G \mathfrak{L}(G) \Leftrightarrow \exp_G^{-1}(Y^*)$ is dense in $\mathfrak{L}(G)$ because the map $\varphi \colon \Delta^* \times \mathfrak{L}(G) \to G$ given by $\varphi(\alpha, r) = \alpha + \exp_G r$ is a local isometry (Proposition 2.14, [1]). However, $\exp_G^{-1}(Y^*)$ is dense in $L(G) \Leftrightarrow Y = \theta_{\Delta^*}^{-1} \exp_G[\exp_G^{-1}(Y^*)]$ is dense in $\mathbb{R}^n \hookrightarrow Y$ has no free summands.

3. We produce a topological isomorphism from a "classical" resolution of a protorus to a new resolution independent of any particular $\Delta \in L(G)$. The group $\frac{\widehat{\Delta}_Y \times \mathfrak{L}(G)}{\Gamma_{L_Y}}$ is LCA by 1 and 2. In the proof of Lemma 4, we found that a subset $U$ of $\widehat{\Delta}_Y$ is open if and only if $U \cap \Delta$ is open in $\Delta$ for all $\Delta \in L_Y$, thus the inclusion map $\Delta^* \times \mathfrak{L}(G) \hookrightarrow \widehat{\Delta}_Y \times \mathfrak{L}(G)$ is a morphism of LCA groups. Define $\eta_Y \colon \frac{\Delta^* \times \mathfrak{L}(G)}{\Gamma_{\Delta^*}} \to \frac{\widehat{\Delta}_Y \times \mathfrak{L}(G)}{\Gamma_{L_Y}}$ to be the morphism induced by inclusion.

We get $\eta[(\alpha, r) + \Gamma_{\Delta^*}] = (\alpha, r) + \Gamma_{L_Y} = 0 \Rightarrow (\alpha, r) \in [\Delta^* \times \mathfrak{L}(G)] \cap \Gamma_{L_Y} \Rightarrow r = -\exp_G^{-1}\alpha$ so $\alpha = -\exp_G r \in \Delta^* \cap \exp_G \mathfrak{L}(G) = \mathbb{Z}_{\Delta^*}$, whence $(\alpha, r) \in \Gamma_{\Delta^*}$; since also $\Gamma_{\Delta^*} \subseteq \Gamma_{L_Y}$ and $\Gamma_{L_Y}$ is closed, it follows that $\eta_Y$ is a well-defined injective morphism.

Next, let $(\beta, s) \in \widehat{\Delta}_Y \times \mathfrak{L}(G)$. By 1, $\beta \in \Delta_F = \overline{F}$ for some free abelian $F$ with $\mathbb{Z}_{\Delta^*} \subseteq F \subseteq Y^* = V_{L_Y}$ and $\Delta_F \cap \exp_G \mathfrak{L}(G) = F$, so that $\beta = \delta + y$ for some $\delta \in \Delta^*$ and $y \in Y^* \subset \exp_G \mathfrak{L}(G)$. Thus, $(y, -\exp_G^{-1} y) \in \Gamma_{L_Y}$ and $(\beta, s) + \Gamma_{L_Y} = (\delta + y, s) - (y, -\exp_G^{-1} y) + \Gamma_{L_Y} = (\delta, s + \exp_G^{-1} y) + \Gamma_{L_Y} = \eta_Y[(\delta, s + \exp_G^{-1} y) + \Gamma_{\Delta^*}]$. This proves that $\eta_Y$ is surjective, whence a topological isomorphism by the open mapping theorem (Theorem 5.29, [5]). Thus, $G \cong_t \frac{\widehat{\Delta}_Y \times \mathfrak{L}(G)}{Y}$ via the diagonal embedding of $Y \cong Y^*$.

4. Let $0 \neq x \in X_\infty$. Then, $X_\infty$ is not $p$-divisible at $x$ if and only if the unique $z_x \in \mathbb{Q}x \cap H$ has $p$-height 0 in $X$. This is equivalent to $\mathbb{Q}x \cap X = \mathbb{Z}_{(p)}z_x$, where $\mathbb{Z}_{(p)}$ is the localization of $\mathbb{Z}$ at $p$. However, $\bigoplus_{z \in H} \frac{\mathbb{Q}z \cap X + \mathbb{Z}^n}{\mathbb{Z}^n} \cong \frac{X}{\mathbb{Z}^n}$, so $\Delta^* \cong \prod_{z \in H} (\frac{\mathbb{Q}z \cap X + \mathbb{Z}^n}{\mathbb{Z}^n})^\vee$ and $\Delta^*$ is $p$-divisible at the point $(\delta_z)_{z \in H}$ with $\delta_z = 0$ for $z \neq z_x$ and $\delta_{z_x} = 1$ where $\overline{\mathbb{Z}(\delta_z)_{z \in H}} \cong (\frac{\mathbb{Q}z_x \cap X + \mathbb{Z}^n}{\mathbb{Z}^n})^\vee$, making appropriate identifications. If $X_\infty$ is $p$-divisible at $0 \neq x \in X_\infty$, then $(\delta_z)_{z \in H}$ has $p$-torsion or $(\overline{\mathbb{Z}(\delta_z)_{z \in H}})_p \cong \widehat{\mathbb{Z}}_p$; but $\theta_{\Delta^*}(z_x) \in X_\infty$ is $p$-divisible so $(\delta_z)_{z \in H}$ is $p$-divisible in $\widehat{\Delta}_{X_\infty}$ according to the topological isomorphism $\eta_{X_\infty}$ in 3. By construction, $X_\infty^*$ is the minimal quotient-divisible torsion-free extension of $X^*$ in $G$, so $\widehat{\Delta}_{X_\infty}$ is the minimal divisible subgroup of $G$ extending $\Delta^*$.

For the last assertion, $\operatorname{tor} G = \operatorname{tor} \widehat{\Delta}_{L(G)}$ by Proposition 5, and this is the same as $\operatorname{tor} \widehat{\Delta}_{\mathbb{Q}^n}$ by 2.

5. $\widehat{\Delta}_{X_\infty} \cong_t \prod_{p \in \mathbb{P}}^{\mathrm{loc}}((\widehat{\Delta}_{X_\infty})_p, \Delta_p^*) \cong_t \varprojlim_{\Delta \in L_{X_\infty}} (\widehat{\Delta}_{X_\infty}/\Delta)$ is an application of (Theorem 3.3, [3]) to the divisible periodic LCA group $\Delta_{X_\infty}$. For the remaining limits, note that in 3 we saw how the lattice $L_Y$ allows us to uniquely determine the topology on $\widehat{\Delta}_Y$ and in turn an associated resolution of $G$. In the range $X \subseteq Y \subseteq \mathbb{Q}^n$, we get a resolution $\frac{\widehat{\Delta}_Y \times \mathfrak{L}(G)}{Y}$ where $Y^*$ is dense in $G$ because $G$ is torus-free. In particular, the topology on the compactly ruled $\widehat{\Delta}_X$ is *coherent* with the collection $L_X$ and it follows $\widehat{\Delta}_X$ is homeomorphic to the the topology on the direct limit of $L_X$ where the upwardly directed partial order of the lattice $L_X$ is preserved. One readily verifies that the topology on the direct limit gives a topological group and our homeomorphism is a topological isomorphism with the tdlc group $\widehat{\Delta}_X$. In parallel, the lattice $M_X \stackrel{\mathrm{def}}{=} \{\mathbb{Z}_\Delta \colon \Delta \in L_X\}$ is isomorphic to the lattice $L_X$, thus the group $X^*$ is isomorphic to the direct limit of $M_X$. The collection $\{\Delta \hookrightarrow G \twoheadrightarrow G/\Delta \colon \Delta \in L_X\}$ of exact sequences dualizes to the collection of exact sequences $\{\mathbb{Z}_\Delta \hookrightarrow X \twoheadrightarrow X/\mathbb{Z}_\Delta \colon \mathbb{Z}_\Delta \in M_X\}$. We conclude that $X \cong \varinjlim_{\Delta \in L_X} (G/\Delta)^\vee$ and, by duality, $X^\vee \cong_t \varprojlim_{\Delta \in L_X}(G/\Delta)$. □

The Structure Theorem for Protori has a number of immediate useful consequences, beginning with the following corollary. Recall the notation $\widehat{\Delta}_p$ introduced for the unique $p$-Sylow subgroup of the $p$-Sylow decomposition of a periodic LCA group $\widehat{\Delta}$ (Theorem 3.3, [3]).

**Corollary 6.** *If $G$ is a torus-free protorus, then $G \cong_t \frac{\widehat{\Delta}_{X_\infty} \times \mathfrak{L}(G)}{X_\infty}$ where $\widehat{\Delta}_{X_\infty} \cong_t \prod_{p \in \mathbb{P}}^{\mathrm{loc}}((\widehat{\Delta}_{X_\infty})_p, \Delta_p)$ is a topological divisible hull of each $\Delta \in L_X$ and $(\widehat{\Delta}_{X_\infty})_p \cong \widehat{\mathbb{Q}}_p^{r_p} \times \mathbb{Z}(p^\infty)^{s_p}$ for some $0 \leq r_p, s_p \in \mathbb{Z}$, $p \in \mathbb{P}$.*

**Proof.** All statements follow directly from Theorem 1 and (Theorem 3.3, Proposition 3.42, [3]). □

For a torus-free protorus $G$, define $\frac{\widehat{\Delta}_{X_\infty} \times \mathfrak{L}(G)}{X_\infty}$ of Corollary 6 to be a *universal resolution* of $G$; the terminology is justified by the fact that $\widehat{\Delta}_{X_\infty}$ is a topological divisible hull of every $\Delta \in L_X$ and previously the only resolutions of protori known to exist were those given in terms of a single element of $L_X$. Recall the notation connotes the closed diagonal embedding $X_\infty \to X_\infty^* \to \Gamma_{L_{X_\infty}} \subseteq \widehat{\Delta}_{X_\infty} \times \mathfrak{L}(G)$. Note that $\frac{\widehat{\Delta}_{L(G)} \times \exp_G \mathfrak{L}(G)}{\widehat{\Delta}_{L(G)} \cap \exp_G \mathfrak{L}(G)}$ is a resolution of $G$ in terms of the canonical subgroups $\widehat{\Delta}_{L(G)}$, generated by all zero-dimensional subgroups (see the proof of Proposition 7), and the path component of 0, $\exp_G \mathfrak{L}(G)$, where $\widehat{\Delta}_{L(G)}$ and $\exp_G \mathfrak{L}(G)$ have their non-locally-compact subspace topologies, and $\widehat{\Delta}_{L(G)} \cap \exp_G \mathfrak{L}(G) \cong \mathbb{Q}^{\dim G}$ (1 in Theorem 1).

**Corollary 7.** *If $G$ is an $n$-dimensional torus-free protorus and $Y$ is torsion-free with $\mathbb{Z}^n \subseteq Y \subseteq \mathbb{Q}^n$, then the map $\varphi_Y \colon \widehat{\Delta}_Y \times \mathfrak{L}(G) \to G$ defined by $\varphi_Y(\alpha, r) = \alpha + \exp_G r$ is a local isometry which is open, continuous, and surjective.*

**Proof.** For each $\Delta \in L(G)$, the map $\varphi_\Delta \colon \Delta \times \mathfrak{L}(G) \to G$ given by $\varphi_\Delta(\alpha, r) = \alpha + \exp_G r$ is a local isometry that is open, continuous, and surjective [1, Proposition 2.2]. By 3 in Theorem 1, the map $\varphi_Y$ has kernel $\Gamma_{L_Y}$ and induces $\frac{\widehat{\Delta}_Y \times \mathfrak{L}(G)}{Y} \cong_t G$, where $Y \to Y^* \to \Gamma_{L_Y} \subseteq \widehat{\Delta}_Y \times \mathfrak{L}(G)$ is the diagonal embedding. However, $L_Y$ is a neighborhood basis at 0 of compact open subgroups for the locally compact topology on $\widehat{\Delta}_Y$ by Proposition 4, and $\varphi_Y$ agrees with $\varphi_\Delta$ for each $\Delta \in L_Y$, so $\varphi_Y$ is a local isometry which is open, continuous, and surjective. □

Define $\widehat{\Delta}_{X_\infty} \times \mathfrak{L}(G)$ to be a *minimal divisible locally compact cover* of $G$. The terminology is justified by the fact that $\widehat{\Delta}_{X_\infty}$ is a topological divisible hull of each $\Delta \in L_X$, the product is locally compact and divisible, and the map $\varphi_{X_\infty} \colon \widehat{\Delta}_{X_\infty} \times \mathfrak{L}(G) \to G$ is open, continuous, and surjective by Corollary 7.

**Corollary 8.** *With the apparatus of Theorem 1 in place for an $n$-dimensional torus-free protorus $G$, set $M_Y \stackrel{\text{def}}{=} \{F \colon F \text{ is a free rank-}n \text{ subgroup of } Y \text{ with } F \subseteq \mathbb{Z}^n \subseteq Y \text{ or } \mathbb{Z}^n \subseteq F \subseteq Y\}$. Then, $L_Y \to M_Y$ given by $\Delta \mapsto \theta_{\Delta^*}^{-1}(\mathbb{Z}_\Delta)$ is bijective and $Y^* = \sum_{\Delta \in L_Y} \mathbb{Z}_\Delta$.*

**Proof.** This follows directly from 1 in Theorem 1, 5 in Theorem 1. □

**Remark 2.** (i) Suppose $G$ is as in Corollary 8 with $Y = \mathbb{Q}^n$. Then, $L_Y = L(G)$ and $M_Y = M_{\mathbb{Q}^n} = \{F \colon F \text{ is a free rank-}n \text{ subgroup of } \mathbb{Q}^n \text{ with } F \subseteq \mathbb{Z}^n \text{ or } \mathbb{Z}^n \subseteq F\}$ with $Y^* = \widehat{\Delta}_{L(G)} \cap \exp_G \mathfrak{L}(G) \cong \mathbb{Q}^n$. (ii) Suppose that $F$ is a free rank-$n$ subgroup of $\widehat{\Delta}_{L(G)} \cap \exp_G \mathfrak{L}(G)$ such that $F \nsubseteq \mathbb{Z}_{\Delta^*}$ and $\mathbb{Z}_{\Delta^*} \nsubseteq F$. Let $\mathbb{Z}^n \subseteq Y \subseteq \mathbb{Q}^n$ with $F \subseteq Y^*$. Then, $F \notin M_Y$ and $\overline{F} \notin L(G)$.

**Proposition 6.** (Protori Lattices) *With the apparatus of Theorem 1 in place for an $n$-dimensional torus-free protorus $G$, $\Delta^* \in L(G)$ with standard representation $\prod_{j=1}^{m} \prod_{p \in \mathbb{P}} \widehat{\mathbb{Z}}(p^{s_p(j)})$, $\mathbb{Z}^n \subseteq Y \subseteq \mathbb{Q}^n$, and $\mathbf{m} \in \mathbb{S}$ with $\mathbf{m}_p = \sup_{z \in H} \{\text{ht}_p^Y(z)\}$, $L_Y$ is isomorphic to the lattice of compact open subgroups of the periodic group $D \stackrel{\text{def}}{=} \prod_{p \in \mathbb{P}}^{\text{loc}} (D_p, C_p)$, where $D_p \stackrel{\text{def}}{=} \prod_{j=1}^{m} D_p(j); C_p \stackrel{\text{def}}{=} \prod_{j=1}^{m} C_p(j);$*

- $D_p(j) = \widehat{\mathbb{Q}}_p$ and $C_p = \widehat{\mathbb{Z}}_p$ if $s_p(j) = \infty$;
- $D_p(j) = \mathbb{Z}(p^\infty)$ and $C_p(j) = p^{-s_p(j)} \widehat{\mathbb{Z}}_p / \widehat{\mathbb{Z}}_p$ if $0 \leqslant s_p(j) < \infty$ and $\mathbf{m}_p = \infty$;
- $D_p(j) = p^{-s_p(j) - \mathbf{m}_p} \widehat{\mathbb{Z}}_p / \widehat{\mathbb{Z}}_p$ and $C_p = p^{-s_p(j)} \widehat{\mathbb{Z}}_p / \widehat{\mathbb{Z}}_p$ if $0 \leqslant s_p(j) + \mathbf{m}_p < \infty$;

$C \stackrel{\text{def}}{=} \prod_{p \in \mathbb{P}} C_p$ *is open in $D$; $D/C \cong \bigoplus_{j=1}^{m} \left\{ \left[ \bigoplus_{s_p(j) + \mathbf{m}_p = \infty} \mathbb{Z}(p^\infty) \right] \oplus \left[ \bigoplus_{s_p(j) + \mathbf{m}_p < \infty} p^{-s_p(j) - \mathbf{m}_p} \mathbb{Z}/\mathbb{Z} \right] \right\}$ is*

*discrete; and the dual of the lattice* $L_Y$ *is isomorphic to the lattice of finite subgroups of* $D/C$.

**Proof.** All periodic LCA groups decompose as a local product as indicated for $D$ by *Braconnier's theorem* (Theorem 3.3, [3]), thus it suffices to determine the $p$-Sylow components $D_p$ and $C_p$. In the proof of Theorem 1, it became evident that the mitigating factor determining the structure of $L_Y$ is the supremum of the $p$-heights in $Y$ of elements in $H$ for each $p \in \mathbb{P}$. $\Delta^* \cong X/\mathbb{Z}^m$ is the base upon which $\widehat{\Delta}_Y$ is formed via the topological isomorphism $\eta_Y$ of 3 in Theorem 1. The last statement is an application of Pontryagin duality to $D$ and its compact open subgroups (Lemma 3.82, [3]). □

**Remark 3.** (i) If $\dim G > 1$, then the lattice of closed subgroups of the protorus $G$, also called the *Chabauty space* $\mathcal{S}(G)$, is equal to the lattice of closed subgroups of $\widehat{\Delta}_{L(G)}$, which is distributive only when $G$ is a solenoid (Theorem 5, [11]). (ii) $\mathcal{S}(\widehat{\Delta}_{L(G)})$ is homeomorphic to $\prod_{\Delta \in L(G)} \mathcal{S}(\Delta)$ (Proposition 1.22, [3]).

## 4. Applications

Protori structure in place, several applications related to morphisms of protori and otherwise follow.

**Lemma 7.** *A morphism* $f_\Delta \colon \Delta_G \to \Delta_H$ *with* $f(\mathbb{Z}_{\Delta_G}) = \mathbb{Z}_{\Delta_H}$ *for some torus-free protori* $G, H$ *and* $\Delta_G \in L(G)$, $\Delta_H \in L(H)$ *extends to an epimorphism* $f \colon G \to H$.

**Proof.** The morphism $\varphi_G \colon \Delta_G \times \mathfrak{L}(G) \to G$ of the *Resolution Theorem* (Proposition 2.2, [1]) is an open map and $\mathbb{Z}_{\Delta_G} \cong_t \exp_G^{-1} \Delta \cong_t \ker \varphi_G$. Let $V \cong_t \mathbb{R}^k$, $0 \leqslant k \in \mathbb{Z}$, denote a real vector space satisfying $\mathfrak{L}(G) = \mathrm{span}_\mathbb{R}(\exp_G^{-1} \Delta) \oplus V$. Then, $G \supseteq \varphi_G(\Delta_G \times V) \cong_t \Delta_G \times V$. The compactness of $G$ implies $k = 0$, so $\exp_G^{-1} \mathbb{Z}_\Delta = \exp_G^{-1} \Delta$ spans $\mathfrak{L}(G)$.

Continuity of $f_\Delta$ with $f(\mathbb{Z}_{\Delta_G}) = \mathbb{Z}_{\Delta_H}$ ensures that $f_\Delta$ is surjective and $\dim_\mathbb{R} \mathfrak{L}(G) = \mathrm{rk}\, \mathbb{Z}_{\Delta_G} \geqslant \mathrm{rk}\, \mathbb{Z}_{\Delta_H} = \dim_\mathbb{R} \mathfrak{L}(H)$. Define $f_\mathbb{R} \colon \mathfrak{L}(G) \to \mathfrak{L}(H)$ by setting $f_\mathbb{R}(\exp_G^{-1}(z)) = \exp_H^{-1}(f(z))$ for $z \in \mathbb{Z}_{\Delta_G}$ and extending $\mathbb{R}$-linearly. Then, $f_\Delta \times f_\mathbb{R} \colon \Delta_G \times \mathfrak{L}(G) \to \Delta_H \times \mathfrak{L}(H)$ is an epimorphism with $(f_\Delta \times f_\mathbb{R})(\Gamma_G) = \Gamma_H$, so $f_\Delta \times f_\mathbb{R}$ induces an epimorphism $\tilde{f} \colon \frac{\Delta_G \times \mathfrak{L}(G)}{\Gamma_G} \to \frac{\Delta_H \times \mathfrak{L}(H)}{\Gamma_H}$ and $\tilde{f}$ in turn induces an epimorphism of protori $f \colon G \to H$ with $f|_{\Delta_G} = f_\Delta$. □

A *projective resolution* of a protorus $G = G_0$ is an exact sequence $K \rightarrowtail P \twoheadrightarrow G$ where $P$ is a torsion-free protorus and $K$ is a torsion-free profinite group. The following result is proven in the narrative immediately following [2, Definitions 8.80].

**Corollary 9.** *A protorus has a projective resolution.*

**Proof.** Let $G$ be a protorus and set $r = \dim G$. By the *Resolution Theorem*, $G$ has a profinite subgroup inducing a torus quotient, which we can take without loss of generality to be $\widehat{\mathbb{Z}}(\vec{\mathbf{n}})$ for some $\vec{\mathbf{n}} \in \mathbb{S}^m$, $m = \mathrm{width}_{nA}\, \widehat{\mathbb{Z}}(\vec{\mathbf{n}})$. Identifying $\mathbb{Z}^r$ in the natural way as a subgroup of $\widehat{\mathbb{Z}}^r$, an isomorphism of free abelian groups $\mathbb{Z}^r \to \mathbb{Z}_{\widehat{\mathbb{Z}}(\vec{\mathbf{n}})}$ extends by continuity to an epimorphism $f_\Delta \colon \widehat{\mathbb{Z}}^r \twoheadrightarrow \widehat{\mathbb{Z}}(\vec{\mathbf{n}})$, thus inducing an exact sequence $K \rightarrowtail \widehat{\mathbb{Z}}^r \twoheadrightarrow \widehat{\mathbb{Z}}(\vec{\mathbf{n}})$ where $K$ is torsion-free profinite. We have $(\widehat{\mathbb{Z}}^r \times \mathbb{R}^r)/\mathrm{diag}(\mathbb{Z}^r) \cong_t P(G) \stackrel{\mathrm{def}}{=} (\mathbb{Q} \otimes G^\vee)^\vee$. By Lemma 7, $f_\Delta$ induces a projective resolution $K \rightarrowtail [\widehat{\mathbb{Z}}^r \times \mathfrak{L}(P(G))]/\Gamma_{P(G)} \twoheadrightarrow [\widehat{\mathbb{Z}}(\vec{\mathbf{n}}) \times \mathfrak{L}(G)]/\Gamma_G$. □

A *completely decomposable group* is a torsion-free abelian group isomorphic to the dual of a completely factorable protorus. An *almost completely decomposable (ACD) group* is a torsion-free abelian group quasi-isomorphic to a completely decomposable group. The Pontryagin dual $G$ of an ACD group is distinguished in the setting of protori by its uniqueness up to topological isomorphism being dependent on a finite factor of an element of $L(G)$.

**Corollary 10.** *If G is a protorus with* $\dim G = \dim_{nA} G$, *then* $G^\vee$ *is an ACD group.*

**Proof.** Let $\Delta_G \in L(G)$. Multiplying $\Delta_G$ by a sufficiently large $N \in \mathbb{Z}$ effects $\text{width}_{nA} N\Delta_G = \dim_{nA} N\Delta_G$. Since $NG = G$, we can assume without loss of generality that $\text{width}_{nA} \Delta_G = \dim_{nA} \Delta_G = \dim G$. Let $\Delta_H$ denote the standard representation for $\Delta_G$ and $\psi \colon \Delta_G \to \Delta_H$ an isomorphism with $\psi(\mathbb{Z}_{\Delta_G}) = \mathbb{Z}_{\Delta_H} \stackrel{\text{def}}{=} \mathbb{Z}e_1 \oplus \cdots \oplus \mathbb{Z}e_{\dim G}$ where $e_j = (0,\ldots,0,1_j,0,\ldots,0)$, $1_j = (1_{jp})_{p\in\mathbb{P}}$, $1 \leqslant j \leqslant m$. By Corollary 5 there is a completely factorable protorus $H$ with $\dim H = \dim G$ and $\Delta_H \in L(H)$. By Lemma 7, there is an epimorphism $\widehat{\psi} \colon G \to H$ extending $\psi$. Symmetrically, there is an epimorphism $\widehat{\eta} \colon H \to G$ extending $\eta \stackrel{\text{def}}{=} \psi^{-1} \colon \Delta_H \to \Delta_G$. It follows that $\widehat{\psi}^\vee \colon H^\vee \to G^\vee$ and $\widehat{\eta}^\vee \colon G^\vee \to H^\vee$ are monomorphisms. By [7, Corollary 6.2.(d)], $G^\vee$ and $H^\vee$ are quasi-isomorphic. It follows that $H^\vee$ is completely decomposable and $G^\vee$ is an ACD group. □

We return to the analysis of morphisms of protori. There is a functor $\mathfrak{L}$ from the category of topological abelian groups to the category of real topological vector spaces [2, Corollary 7.37]: for a morphism $f \colon G \to H$ of topological abelian groups, the map $\mathfrak{L}(f) \colon \mathfrak{L}(G) \to \mathfrak{L}(H)$ given by $\mathfrak{L}(f)(r) = f \circ r$ is a morphism of real topological vector spaces satisfying $\exp_H \circ \mathfrak{L}(f) = f \circ \exp_G$.

**Proposition 7.** *A morphism* $G \to H$ *between torus-free protori restricts to morphisms* $\widehat{\Delta}_{L(G)} \to \widehat{\Delta}_{L(H)}$, $\exp_G \mathfrak{L}(G) \to \exp_H \mathfrak{L}(H)$, *and a continuous map* $X_G \to X_H$.

**Proof.** Let $D$ be a profinite subgroup of $G$. If $\Delta \in L(G)$, then $\Delta + D$ is profinite because it is compact and zero-dimensional: the addition map $\Delta \times D \twoheadrightarrow \Delta + D$ is a continuous epimorphism and the kernel $K$ is closed (whence profinite), so we get an exact sequence $K \rightarrowtail \Delta \times D \twoheadrightarrow \Delta + D$, whence $\dim(\Delta + D) = \dim(\Delta \times D) - \dim K = \dim \Delta + \dim D - \dim K = 0$ by Lemma 1. The natural map $G/\Delta \to G/(\Delta + D)$ is surjective, so $\Delta + D \in L(G)$. Hence, $D \subseteq \Delta + D \subseteq \widehat{\Delta}_{L(G)}$. We conclude that $\widehat{\Delta}_{L(G)} = \sum\{D \colon D$ a profinite subgroup of $G\}$, and similarly for $\widehat{\Delta}_{L(H)}$. In particular, $\widehat{\Delta}_{L(G)}$ *contains all profinite subgroups of G*; similarly for $\widehat{\Delta}_{L(H)}$.

Let $f$ denote a morphism $G \to H$. If $\Delta \in L(G)$, then $K = \ker f \cap \Delta$ is profinite, so $\Delta/K \cong_t f(\Delta)$ is profinite. Thus, $f(\Delta) \subseteq \widehat{\Delta}_{L(H)}$. It follows that $f(\widehat{\Delta}_{L(G)}) \subseteq \widehat{\Delta}_{L(H)}$. In addition, $\exp_H \circ \mathfrak{L}(f) = f \circ \exp_G$ implies that $f[\exp_G \mathfrak{L}(G)] \subseteq \exp_H \mathfrak{L}(H)$. Lastly, Theorem 1 gives that $f(X_G) = f(\widehat{\Delta}_{L(G)} \cap \exp_G \mathfrak{L}(G)) \subseteq f(\widehat{\Delta}_{L(G)}) \cap f(\exp_G \mathfrak{L}(G)) \subseteq \widehat{\Delta}_{L(H)} \cap \exp_H \mathfrak{L}(H) = X_H$. □

**Remark 4.** *The continuous map* $X_G \to X_H$ *in Proposition 7 is not, in general, a homomorphism of torsion-free abelian groups.*

**Proposition 8.** *For a morphism* $f \colon G \to H$ *of torus-free protori there exist* $\Delta_G \in L(G)$, $\Delta_H \in L(H)$ *such that* $f$ *lifts to a product map* $f|_{\Delta_G} \times \mathfrak{L}(f) \colon \Delta_G \times \mathfrak{L}(G) \to \Delta_H \times \mathfrak{L}(H)$.

**Proof.** Let $\Delta_G \in L(G)$. By Proposition 7, $f(\widehat{\Delta}_{L(G)}) \subseteq \widehat{\Delta}_{L(H)}$. By Theorem 1, $\widehat{\Delta}_{L(H)} = \bigcup_{\Delta \in L(H)} \Delta$. Each $\Delta \in L(H)$ is open in $\widehat{\Delta}_{L(H)}$ because the intersection of any two elements of $L(H)$ is an element of $L(H)$ with finite index in any other element of $L(H)$ containing it (Proposition 2.1.2, [8]). By Proposition 7, $f(\Delta_G) \subseteq \widehat{\Delta}_{L(H)}$. Because $f(\Delta_G)$ is compact and the elements of $L(H)$ are open in $\widehat{\Delta}_{L(H)}$, there are finitely many elements of $L(H)$ which cover $f(\Delta_G)$; let $\Delta_H \in L(H)$ denote the sum of these elements. Then, $f(\Delta_G) \subseteq \Delta_H$. Since $\exp_H \circ \mathfrak{L}(f) = f \circ \exp_G$, it follows that $f|_{\Delta_G} \times \mathfrak{L}(f) \colon \Delta_G \times \mathfrak{L}(G) \to \Delta_H \times \mathfrak{L}(H)$ is a lifting of $f \colon G \to H$. □

A morphism of torus-free protori lifts to one between the minimal divisible locally compact covers:

**Theorem 2.** (Structure Theorem for Morphisms) *A morphism* $f \colon G \to H$ *of torus-free protori with duals $X$ and $Y$ lifts to a product map* $f|_{\widehat{\Delta}} \times f_\mathfrak{L} \colon \widehat{\Delta}_{X_\infty} \times \mathfrak{L}(G) \to \widehat{\Delta}_{Y_\infty} \times \mathfrak{L}(H)$.

**Proof.** This follows from Proposition 8 because $\widehat{\Delta}_{X_\infty} = \sum_{\Delta \in L_{X_\infty}} \Delta$. □

**Funding:** This research received no external funding.

**Acknowledgments:** First and foremost, I am grateful to my mentor, Adolf Mader, who came out of retirement to support my reintroduction to mathematics research 25 years after serving as my Ph.D. adviser; the content of this paper and the opportunity to pursue my passion of doing mathematics simply would not exist if not for him. I would like to thank Keith Conrad of the University of Connecticut for his generous help answering all my number theory questions, graciously taking the time to read and respond to my meandering emails. Lastly, thanks are due to Karl Hofmann for his encouragement as I learned the fundamentals from the masterful treatise he coauthored with Sidney Morris: *The Structure of Compact Groups*.

**Conflicts of Interest:** The author declares no conflict of interest.